\documentclass{amsart}
\usepackage{enumerate,amssymb,amsmath,mathrsfs}
\usepackage[dvips]{graphics}
\usepackage[dvips]{graphicx}
\usepackage{epsfig,enumerate,mathrsfs}
\usepackage{color}

\numberwithin{equation}{section}
\newtheorem{theorem}{Theorem}[section]
\newtheorem{proposition}[theorem]{Proposition}
\newtheorem{lemma}[theorem]{Lemma}
\newtheorem{corollary}[theorem]{Corollary}
\newtheorem{remark}[theorem]{Remark}
\newtheorem{example}{Example}

\DeclareMathOperator{\Tr}{Tr}

\newcommand{\Ss}{\mathcal{S}}

\newcommand{\la}{\lambda}
\newcommand{\Aa}{\mathcal{A}}
\newcommand{\Bb}{\mathcal{B}}

\newcommand{\z}{\zeta}

\newcommand{\Tt}{\mathcal{T}}
\newcommand{\Ll}{\mathcal{L}}

\newcommand{\R}{\mathbb R}
\newcommand{\Z}{\mathbb Z}
\newcommand{\Cc}{\mathcal{C}}

\def\N{\mathbb{N}}
\def\C{\mathbb{C}}
 
\def\Id{\mathrm{Id}}
\def\dom{\mathrm{dom}}

\def\KD3{\mathrm{KD}^3}

\def\pP{\mathscr{P}}
\def\lL{\mathscr{L}}

\def\max{\mathrm{max}}
\def\Ei{\mathrm{Ei}}

\begin{document}

\title[Unusual properties of zeta functions on cones]
{Exotic expansions and pathological properties of
$\zeta$-functions on conic manifolds}





\author{Klaus Kirsten}
\address{Department of Mathematics\\ Baylor University\\
         Waco\\ TX 76798\\ U.S.A. }
\email{Klaus$\_$Kirsten@baylor.edu}

\author{Paul Loya}
\address{Department of Mathematics \\Binghamton University\\
Binghamton\\NY 13902\\ U.S.A. }
\email{paul@math.binghamton.edu}

\author{Jinsung Park}
\address{School of Mathematics\\ Korea Institute for Advanced Study\\
207-43\\ Cheongnyangni 2-dong\\ Dongdaemun-gu\\ Seoul 130-722\\
Korea } \email{jinsung@kias.re.kr}




\begin{abstract}
{We give a complete classification and present new exotic
phenomena of the meromorphic structure of $\z$-functions
associated to general self-adjoint extensions of Laplace-type
operators over conic manifolds. We show that the meromorphic
extensions of these $\z$-functions have, in general, countably
many logarithmic branch cuts on the nonpositive real axis and
unusual locations of poles with arbitrarily large multiplicity.
The corresponding heat kernel and resolvent trace expansions also
exhibit exotic behaviors with logarithmic terms of arbitrary
positive and negative multiplicity. We also give a precise
algebraic-combinatorial formula to compute the coefficients of the
leading order terms of the singularities.}
\end{abstract}

\maketitle


\section{Introduction}


In this paper we give a complete classification of the meromorphic
structure of $\zeta$-functions associated to conic manifolds; that
is, general self-adjoint extensions of Laplace-type operators on
conic manifolds introduced by Cheeger \cite{ChJ79,Ch2}. In
particular, we prove that such $\zeta$-functions exhibit
pathological meromorphic properties. Before giving a synopsis of
these pathological properties, recall that the $\zeta$-function
$\zeta(s,\Delta)$ of a Laplacian $\Delta$ over a smooth closed
manifold has a meromorphic extension to all of $\C$ with only
simple poles at $s = \frac{n-k}{2}\notin -\N_0$ with $n$ the
dimension of the manifold and $k \in \N_0:=\{0,1,2,\ldots \}$
\cite{MiS-PlA49,MinS49,SeR67,WeyH49}. The situation is completely
different for conic manifolds. We show that the $\zeta$-function
associated to a general self-adjoint extension of a Laplace-type
operator on a conic manifold has, as a general rule (except for
very special cases, e.g.\ the Friedrichs extension), in addition
to the singularities at $s = \frac{n-k}{2} \notin -\N_0$ for $k
\in \N_0$, the following properties:
\begin{enumerate}[\it (1)]
\item It can have countably many poles of \emph{arbitrarily} high
multiplicity at ``unusual" locations on the negative real axis;
that is, at points not of the form $s = \frac{n-k}{2}$.


\item It can have \emph{countably} many logarithmic singularities
at ``unusual" locations.

\item The singularities in {\it (1)} and {\it (2)} can occur for
the \emph{same} $\zeta$-function and at the \emph{same} ``unusual"
locations. Moreover, we also give an elementary and explicit
algebraic-combinatorial recipe to compute the exact locations and
leading coefficients of the ``unusual" poles and logarithmic
singularities.
\end{enumerate}


In fact, the explicit computation of these exotic singularities is
so straightforward (see Section \ref{ssec-ex}) that for low
dimensions we can find the structure of zeta functions quickly. We
also remark that one can always conjure up \emph{artificial} zeta
functions having {\it (1)} and {\it (2)}, but for \emph{natural}
(geometric) zeta functions, properties {\it (1)} and {\it (2)}
seem to have no parallels in the differential geometry literature.

\subsection{A simple example} \label{ssec-egR2}

Here is a surprising, and completely natural, example of a
$\zeta$-function which has no meromorphic extension to all of
$\C$. We first review conic manifolds. Let $M$ be an
$n$-dimensional compact manifold with boundary $\Gamma$ and let
$g$ be a smooth Riemannian metric on $M \setminus \Gamma$. We
assume that near $\Gamma$ there is a collared neighborhood
\[
\mathcal{U} \cong [0,\varepsilon)_r \times \Gamma ,
\]
where $\varepsilon > 0$ and the metric $g$ is of product type
$dr^2 + r^2 h$ with $h$ a metric over $\Gamma$. Such a metric is
called a \emph{conic metric} and $M$ is called a \emph{conic
manifold}, ideas introduced by Cheeger \cite{ChJ79,Ch2} (cf.\
\cite{MM}). Using a Liouville transformation over the collar
$\mathcal{U}$ as in \cite{BS1}, we can identify $L^2(M,dg)$ with
$L^2(M, dr dh)$ and the scalar Laplacian $\Delta_g$ can be
identified with
\begin{equation} \label{Deltaconic}
\Delta_g \big|_{\mathcal{U}} = - \partial_r^2 +
\frac{1}{r^2}A_{\Gamma},\qquad \text{with}\quad A_{\Gamma} =
\Delta_\Gamma + \Big(\frac{1 - n}{2}\Big)
              \Big(1 + \frac{1 - n}{2}\Big) ,
\end{equation}
where $\Delta_\Gamma$ is the Laplacian over $\Gamma$. Notice that
$A_\Gamma \geq -\frac14$ because the function $x ( 1 + x)$ has the
minimum value $- \frac14$ (when $x = -\frac12$). Let us now assume
that $n = 2$ so that
\[
\Delta_g \big|_{\mathcal{U}} = - \partial_r^2 +
\frac{1}{r^2}A_{\Gamma},\qquad \text{with}\quad A_{\Gamma} =
\Delta_\Gamma - \frac{1}{4} .
\]
We remark that the term $- \frac{1}{4 r^2}$ can be considered a
``singular potential," and such Laplacians and their self-adjoint
extensions have been studied by physicists since the 70's
\cite{BuW-GeF85,CoS-HoB02,FrW-LaD-SpR71,RadC75}. Note that
$\Delta_\Gamma$ always has the eigenvalue $0$. Let us assume that
$0$ is the only eigenvalue (counting multiplicity) of
$\Delta_\Gamma$ in the interval $[0,1)$; this is the case for the
Euclidean Laplacian on a punctured region in $\R^2$, for in this
case, $\Delta_\Gamma$ is just the Laplacian on the unit circle
which has eigenvalues $\{k^2\, |\, k \in \Z\}$. Then $A_{\Gamma}$
has exactly one eigenvalue in the interval $[-\frac14, \frac34)$,
the eigenvalue $- \frac14$, and (see Section \ref{s:SAE})
$\Delta_g$ has many different self-adjoint extensions, each of
which is parameterized by an angle $\theta \in [0,\pi)$ (cf.\
\cite{KocA79,KaB-StU91,GM}). It turns out that $\theta =
\frac{\pi}{2}$ corresponds to the so-called Friedrichs extension
\cite{BS1}, and any extension has a discrete spectrum
\cite{BLeM97}. Consider any one of the extensions, say
$\Delta_\theta$ with $\theta \in [0,\pi)$, and form the
corresponding $\zeta$-function
\[
\zeta(s ,\Delta_\theta) := \sum_{\mu_j \ne 0} \frac{1}{\mu_j^s},
\]
where the $\mu_j$'s are the eigenvalues of $\Delta_\theta$. The
surprising fact is that the meromorphic extension of \emph{every}
such $\zeta$-function corresponding to an angle $\theta \in
[0,\pi)$ \emph{except} $\theta = \frac{\pi}{2}$ (the Friedrich's
extension) has a logarithmic singularity at $s = 0$. More
precisely, as a consequence of Theorem \ref{main thm2}, for
$\theta \ne \frac{\pi}{2}$, we can write
\begin{equation}
\z(s,\Delta_\theta) = \z_{\mathrm{reg}}(s,\Delta_\theta) - \frac{
\sin (\pi s)}{\pi}e^{-2 s \kappa} \log s \,\, ,\label{add1}
\end{equation}
where $\kappa = \log 2 - \gamma - \tan \theta$ with $\gamma$ the
Euler-Mascheroni constant and $\z_{\mathrm{reg}}(s,\Delta_\theta)$
has a meromorphic extension over $\C$ with the ``regular" simple
poles at the well-known values $s = \frac12 - k$ for $k \in
\mathbb{N}_0$.

\subsection{Operators on conic manifolds} \label{sec-RSO}
Br\"uning and Seeley's regular singular operators \cite{BS1}
generalize the example \eqref{Deltaconic} of the Laplacian on a
conic manifold as follows. Let $M$ be an $n$-dimensional compact
manifold with boundary $\Gamma$ and we assume that near $\Gamma$
there is a collared neighborhood $\mathcal{U}$ such that
\[
\mathcal{U} \cong [0,\varepsilon)_r \times \Gamma ,
\]
where $\varepsilon > 0$ and the metric of $M$ is of product type
$dr^2 + h$ with $h$ a metric over $\Gamma$. Let $E$ be a Hermitian
vector bundle over $M$ and let $\Delta$ be a second order
\emph{regular singular operator} acting on $C^\infty_c(M\setminus
\Gamma, E)$; this means that $\Delta$ is an elliptic symmetric
nonnegative second order differential operator such that the
restriction of $\Delta$ to $\mathcal{U}$ has the form
\begin{equation} \label{D}
- \partial_r^2 + \frac{1}{r^2} A_\Gamma,
\end{equation}
where $A_\Gamma$ is a Laplace-type operator over $\Gamma$ with
$A_\Gamma \geq -\frac14$. (The condition $A_\Gamma \geq -\frac14$
is necessary otherwise $\Delta$ is not bounded below
\cite{BS1,CaC83}.) Laplacians on forms and squares of Dirac
operators on conic manifolds
\cite{BS1,ChJ79,Ch2,Chou,BKirK01,ULeMo,LMP} are examples of second
order regular singular operators. We can also deal with the case
when $M$ has boundary components up to which $\Delta$ is smooth;
at such components, put local boundary conditions.

As shown by Mooers in \cite[Prop.\ 2.4]{MoE99} (cf.\ Gil and
Mendoza \cite{GM}), the self-adjoint extensions of $\Delta$ are
parameterized by Lagrangian subspaces in the eigenspaces of
$A_\Gamma$ with eigenvalues in the interval $[- \frac14,
\frac34)$. To describe these extensions, denote by
\[
-\frac14 = \underbrace{\lambda_1 = \lambda_2 = \cdots =
\lambda_{q_0}}_{ = - \frac14} < \underbrace{\lambda_{q_0+1} \leq
\lambda_{q_0 + 2} \leq \cdots \leq \la_{q_0 + q_1}}_{- \frac14 <
\la_\ell < \frac34}
\]
the spectrum of $A_\Gamma$ in $[- \frac14, \frac34)$ and by
$\{\phi_\ell\}$ the associated eigenvectors, and define
\[
V \cong \bigoplus_{- \frac14 \leq \la_\ell < \frac34} E_\ell
\oplus E_\ell \cong \C^{2 q} \qquad \text{with}\ \ E_\ell :=
\mathrm{span} \left\{ \phi_\ell \right\}\ \ \text{and}\ \ q = q_0
+ q_1;
\]
see Section \ref{s:SAE} for a more precise description of $V$. We
can endow $V$ with a symplectic structure as described in Section
\ref{s:SAE}. Then the self-adjoint extensions of $\Delta$ are in a
one-to-one correspondence to the Lagrangian subspaces in $V$.
Given a Lagrangian subspace $L$ in $V$, we denote by $\Delta_L$
the self-adjoint extension corresponding to $L$.

One of the natural questions for a given self-adjoint extension
$\Delta_L$ is whether the $\z$-function of $\Delta_L$,
$\z(s,\Delta_L)$, would have a meromorphic extension over $\C$ and
if so, what the pole structure is. Here the $\z$-function of
$\Delta_L$ is defined by
\begin{equation}\label{e:def-zeta}
\z(s,\Delta_L) = \sum_{\mu_j \neq 0} \frac{1}{\mu_j^s}
\end{equation}
for $\Re s\gg 0$ where $\mu_j$'s are the eigenvalues of
$\Delta_L$. The $\z$-function has been studied in many papers for
the Friedrich's extension or with the condition $A_\Gamma \geq
\frac34$, which implies that $\Delta$ is essentially self-adjoint
\cite{BKD,BGKE,BS1,CaC83,CaC88,Ch1,Ch2,CoG-ZeS97,G,GL,Les,LoPR02,LMP,SprM05}.
For scalable extensions (extensions for which the domain is
invariant under $r \mapsto c r$), the $\zeta$-function has been
studied in \cite{BLeM97}. There are no ``unusual" phenomena with
these cases. Finally, for arbitrary self-adjoint extensions with
$A_\Gamma \geq - \frac14$, the $\zeta$-function has been studied
by Falomir, Muschietti and Pisani\ \cite{FMP} (see also
\cite{FMPS,FPW}) for one-dimensional Laplace-type operators over
$[0,1]$ and by Mooers \cite{MoE99} who was the first to study the
general case of operators over manifolds and who noticed the
presence of ``unusual" poles. However, the works \cite{FMP,MoE99}
only imply the existence of \emph{simple} ``unusual" poles and do
not imply the existence of poles of arbitrary order nor of
logarithmic singularities of the $\zeta$-function.

We now outline this paper. We begin in Section \ref{s:statement}
by giving the statement of our main result, Theorem
\ref{thm-main}, and we also illustrate the ease of applying the
main result by giving examples; in particular, we re-derive the
main result of \cite{FMP} and we show that poles of arbitrary
order and countably many logarithmic singularities show up even
for simple situations. We also show how our theorem simplifies
when we make assumptions on the self-adjoint extensions and we
present corresponding resolvent and heat kernel expansions. Mooers
\cite{MoE99} description of self-adjoint extensions as Lagrangian
subspaces plays a key r\^ole in the proof of our main result and
because of this reason, in Section \ref{s:SAE} we briefly review
this important topic. The main technical task of the proof of our
main theorem is the explicit form of the parametrix of the
resolvent of $\Delta_L$ near the boundary. This is handled by
solving model problems explicitly over a finite interval employing
Theorem \ref{thm-res}, which is a new representation of the
resolvent in terms of an implicit eigenvalue equation, and the
contour integration method \cite{BKirK01,KM1,KM2}. The
presentation and details of the solutions of the model problems
are given in Sections \ref{s:modelproblems}, \ref{s:modelzeta} and
\ref{s:model}. The results for the model problems and a parametix
construction then enable us to prove all the theorems listed in
Section \ref{s:statement} below. This is described in Section
\ref{s:meromor}.

\section{Statement and examples of results}
\label{s:statement}

\subsection{Statement of main result}

Fix a Lagrangian $L \subset V$ and hence a self-adjoint extension
$\Delta_L$ of $\Delta$ (we use the notation from Section
\ref{sec-RSO}). In Section \ref{s:SAE}  we show that $L$ can be
described by $q \times q$ matrices $\Aa$ and $\Bb$ having the
property that the rank of the $q \times 2 q$ matrix
$\begin{pmatrix} \Aa & \Bb
\end{pmatrix}$ is $q$ and $\Aa' \, \Bb^*$ is self-adjoint where
$\Aa'$ is the matrix $\Aa$ with the first $q_0$ columns multiplied
by $-1$ (conversely, any such $\Aa$ and $\Bb$ define a
Lagrangian). Before stating the main result which describes the
exact structure of $\zeta(s,\Delta_L)$, we apply a straightforward
three-step algorithm to $\Aa$ and $\Bb$ that we need for the
statement.

{\bf Step 1:} First, we define the function
\begin{equation} \label{polyp}
p(x,y) := \det
\begin{pmatrix} \Aa & \Bb\\
\begin{array}{cccc}
x\, \Id_{q_0} & 0 & 0 & 0\\
0 & \tau_1\, y^{2 \nu_1} & 0 & 0\\
0 & 0 & \ddots & 0\\
0 & 0 & 0 & \tau_{q_1}\, y^{2 \nu_{q_1}}
\end{array}
&  \Id_q &
\end{pmatrix},
\end{equation}
where $\Id_k$ denotes the $k \times k$ identity matrix and where
\[
\nu_j := \sqrt{\la_{q_0 + j} + \frac14}\ , \quad \tau_j =
2^{2\nu_j}\frac{ \Gamma(1 + \nu_j)}{ \Gamma(1 - \nu_j)},\qquad j =
1,\ldots, q_1.
\]
For specific $\Aa$ and $\Bb$, $p(x,y)$ is explicitly computable
``by hand"; we shall give some examples in Subsection
\ref{ssec-ex}. Expanding the determinant using one's favorite
method, we can write $p(x,y)$ as a ``polynomial"
\[
p(x,y) = \sum a_{j \alpha} \, x^j\, y^{2 \alpha} ,
\]
where the $\alpha$'s are linear combinations of $\nu_1, \ldots,
\nu_{q_1}$ and the $a_{j \alpha}$'s are constants. Let $\alpha_0$
be the smallest of all $\alpha$'s with $a_{j \alpha} \ne 0$ and
let $j_0$ be the smallest of all $j$'s amongst the $a_{j \alpha_0}
\ne 0$. Then factoring out the term $a_{j_0 \alpha_0}\, x^{j_0}\,
y^{2 \alpha_0}$ in $p(x,y)$ we can write $p(x,y)$ in the form
\begin{equation} \label{pxy}
p(x,y) = a_{j_0 \alpha_0}\, x^{j_0}\, y^{2 \alpha_0} \Big( 1 +
\sum  b_{k \beta} \, x^k \, y^{2 \beta} \Big)
\end{equation}
for some constants $b_{k \beta}$ (equal to $a_{k \alpha}/a_{j_0
\alpha_0}$).

{\bf Step 2:} Second, using the power series $\log (1 + z) =
\sum_{k = 1}^\infty \frac{(-1)^{k-1}}{k} z^k$ with $z = \sum  b_{k
\beta} \, x^k \, y^{2 \beta}$ for a sufficiently small $|z|$, we
can write
\begin{equation} \label{clx}
\log \Big( 1 + \sum  b_{k \beta} x^k y^{2 \beta} \Big) = \sum
c_{\ell \xi} \, x^\ell \, y^{2 \xi}
\end{equation}
for some constants $c_{\ell \xi}$. We emphasize that for specific
$\Aa$ and $\Bb$, all the coefficients $c_{\ell \xi}$ are
explicitly computable ``by hand" when $q$ is small (see the
examples in Subsection \ref{ssec-ex}) and easily with a computer
for $q$ large. With a little thought, one can see that the $\xi$'s
appearing in \eqref{clx} are nonnegative, countable and approach
$+\infty$ unless $\beta=0$ is the only $\beta$ occurring in
\eqref{clx}, in which case $\xi=0$ occurs in \eqref{clx}. Also the
$\ell$'s with $c_{\ell \xi} \ne 0$ for a fixed $\xi$ are bounded
below.

{\bf Step 3:} Third, for each $\xi$ appearing in \eqref{clx},
define
\begin{equation} \label{ellxi}
p_\xi := \min \{\ell \leq 0 \, |\, c_{\ell \xi} \ne 0\}\qquad
\text{and}\qquad \ell_\xi := \min \{\ell > 0 \, |\, c_{\ell \xi}
\ne 0\},
\end{equation}
when these numbers are actually defined, that is, whenever the
sets $\{\ell \leq 0 \, |\, c_{\ell \xi} \ne 0\}$ and $\{\ell
> 0 \, |\, c_{\ell \xi} \ne 0\}$, respectively, are nonempty. We
now define
\begin{equation} \label{Spl}
{\pP} := \{ \xi \, |\, p_\xi\ \text{is defined}\} \qquad
\text{and}\qquad {\lL} := \{\xi \, |\, \ell_\xi\ \text{is
defined}\}.
\end{equation}
The following theorem is our main result.

\begin{theorem} \label{thm-main}
Let $L \subset V$ be an \textbf{arbitrary} Lagrangian subspace of
$V$ and define $\mathscr{\pP}$ and $\mathscr{\lL}$ as in
\eqref{Spl} from the matrices $\Aa$ and $\Bb$ defining $L$. Then
the $\zeta$-function $\zeta(s,\Delta_L)$ extends from $\Re s
> \frac{n}{2}$ to a holomorphic function on $\C \setminus
(-\infty,0]$. Moreover, $\zeta(s,\Delta_L)$ can be written in the
form
\[
\zeta(s,\Delta_L) = \zeta_{\mathrm{reg}}(s,\Delta_L) +
\zeta_{\mathrm{sing}}(s,\Delta_L),
\]
where $\zeta_{\mathrm{reg}}(s,\Delta_L)$ has possible ``regular"
poles at the ``usual" locations $s = \frac{n - k}{2}\notin -\N_0$
for $k \in \N_0$ and at $s = 0$ if $\dim \Gamma
> 0$, and where $\zeta_{\mathrm{sing}}(s,\Delta_L)$ has the
following expansion:
\begin{multline} \label{zetasing}
\zeta_{\mathrm{sing}}(s,\Delta_L) = \frac{\sin (\pi s)}{\pi}
\bigg\{ (j_0 - q_0) e^{-2 s (\log 2 - \gamma)} \log s + \sum_{\xi
\in
\mathscr{\pP}} \frac{f_\xi(s)}{(s + \xi)^{|p_\xi| + 1}}\\
+ \sum_{\xi \in \mathscr{\lL}} g_\xi(s) \log (s + \xi) \bigg\},
\end{multline}
where $j_0$ appears in \eqref{pxy} and $f_\xi(s)$ and $g_\xi(s)$
are entire functions of $s$ such that
\[
f_\xi(-\xi) = (-1)^{|p_\xi|+1} c_{p_\xi \xi} \, \xi \,
\frac{|p_\xi|!}{2^{|p_\xi|}} \,
\]
and near $s = - \xi$,
\[
g_\xi(s) = \begin{cases} c_{\ell_0, 0} \, \frac{2^{\ell_0}}{
(\ell_0 - 1)!} s^{\ell_0} + \mathcal{O}(s^{\ell_0 + 1}) & \text{if
\ \  $\xi = 0$,}\\ - c_{\ell_\xi \xi} \, \frac{\xi
2^{\ell_\xi}}{(\ell_\xi - 1)!} (s + \xi)^{\ell_\xi - 1} +
\mathcal{O}((s + \xi)^{\ell_\xi}) &  \text{if \ \ $\xi > 0$.}
\end{cases}
\]
\end{theorem}

\begin{remark} \label{rmk-main} \em
The zeta function $\zeta_{\mathrm{reg}}(s,\Delta_L)$ will only
have possible poles at $s = \frac{n-k}{2}\notin -\N_0$ in the case
that $\Gamma$ is the only boundary component of $M$ and the
residue of $\zeta_{\mathrm{reg}}(s,\Delta_L)$ at $s = 0$ is given
by $\mathrm{Res}_{s = 0} \zeta_{\mathrm{reg}}(s,\Delta_L) = -
\frac12 \mathrm{Res}_{s = - \frac12} \zeta (s, A_\Gamma )$; in
particular, this vanishes if $\zeta (s, A_\Gamma )$ is in fact
analytic at $s = -\frac12$. Later in Theorems
\ref{thm-restracegen} and \ref{t:traceM}, we shall present
corresponding resolvent and heat kernel expansions. In general,
$\mathscr{\lL}$ may contain the origin, so that the logarithm
part, which has the branch cut at $s=0$, is given by
\[
\frac{\sin (\pi s)}{\pi} \big( (j_0 - q_0) e^{-2 s (\log 2 -
\gamma)}+ g_0(s)\big) \log s
\]
where $g_0(s)$ depends on $\Aa,\Bb$.  Also, it is easy to check
that when there are no $-\frac14$ eigenvalues, then there are no
logarithmic singularities and the ``unusual" poles occur with
multiplicity at most one. Finally, the expansion \eqref{zetasing}
means that for any $N \in \N$,
\begin{multline*}
\zeta_{\mathrm{sing}}(s,\Delta_L) = \frac{\sin (\pi s)}{\pi}
\bigg\{ (j_0 - q_0) e^{-2 s (\log 2 - \gamma)} \log s  + \sum_{\xi
\in \mathscr{\pP} ,\, \xi \leq N} \frac{f_\xi(s)}{(s +
\xi)^{|p_\xi| + 1}}\\ + \sum_{\xi \in \mathscr{\lL} ,\, \xi \leq
N} g_\xi(s) \log (s + \xi) \bigg\} + F_N(s),
\end{multline*}
where $F_N(s)$ is holomorphic for $\Re s \geq - N$.
\end{remark}

\subsection{Examples} \label{ssec-ex}

Via examples we show the ease and efficiency at which Theorem
\ref{thm-main} computes the exact meromorphic structure of
$\zeta_{\mathrm{sing}}(s,\Delta_L)$ (note that
$\zeta_{\mathrm{reg}}(s,\Delta_L)$ is ``uninteresting," which is
why we focus on $\zeta_{\mathrm{sing}}(s,\Delta_L)$).

\begin{example} \label{eg-FMP} \em (Taken from \cite{FMP}.)
The paper by Falomir \emph{et al.}  \cite{FMP} (along with Mooers'
\cite{MoE99}) is in many ways the inspiration for our paper and is
the very first paper to find explicit formulas for the ``unusual"
poles of Laplacians; cf.\ \cite{FPW} for the infinite interval and
\cite{FMPS} for squares of $2 \times 2$ systems. \cite{FMP}
studies the operator
\[
\Delta = - \frac{d^2}{dr^2} + \frac{1}{r^2}\, \la\qquad
\text{over}\ \ [0,1]
\]
with $- \frac14 \leq \la < \frac34$. In this case, $V = \C^2$,
therefore Lagrangians $L \subset \C^2$ are determined by $1 \times
1$ matrices (numbers) $\Aa = \alpha$ and $\Bb = \beta$. Fix such
an $(\alpha,\beta)$; we shall determine the strange singularity
structure of $\zeta(s,\Delta_L)$. Let us assume that $- \frac14 <
\la < \frac34$ so there are no $- \frac14$ eigenvalues (we will
come back to the $\la = - \frac14$ in a moment). Then with $\nu :=
\sqrt{\la + \frac14}$ and $\tau := 2^{2\nu}\frac{ \Gamma(1 +
\nu)}{ \Gamma(1 - \nu)}$,
\[
p(x,y) := \det
\begin{pmatrix} \alpha & \beta \\
\tau \, y^{2 \nu} & 1 \end{pmatrix} = \alpha - \beta \, \tau\,
y^{2 \nu} = \alpha\Big( 1 - \frac{\tau \beta}{\alpha} y^{2 \nu}
\Big),
\]
where we assume that $\alpha,\beta \ne 0$ (the $\alpha = 0$ or
$\beta = 0$ cases can be handled easily), and we write $p(x,y)$ as
in \eqref{pxy}. Forming the power series \eqref{clx}, we see that
\[
\log \Big( 1 - \frac{\tau \beta}{\alpha} y^{2 \nu} \Big) = \sum_{k
= 1}^\infty \frac{(-1)^{k-1}}{k} \Big(- \frac{\tau \beta}{\alpha}
y^{2 \nu} \Big)^k = \sum_{k = 1}^\infty c_{0, \nu k}\, x^0 y^{2
\nu k},
\]
where $c_{0, \nu k} = - \frac{1}{k} \Big( \frac{\tau
\beta}{\alpha} \Big)^k$. Using the definitions \eqref{ellxi} and
\eqref{Spl} for $p_\xi$, $\ell_{\xi}$, $\mathscr{\pP}$, and
$\mathscr{\lL}$, we immediately see that $\ell_{\nu k}$ never
exists so $\mathscr{\lL} = \varnothing$, while
\[
p_{\nu k} = \min \{\ell \leq 0\, |\, c_{\ell , \nu k} \ne 0\} = 0
\ , \quad \mathscr{\pP} = \{\nu \, k\, |\, k \in \N\}.
\]
Therefore, by Theorem \ref{thm-main},
\[
\zeta_{\mathrm{sing}}(s,\Delta_L) = \frac{\sin (\pi s)}{\pi}
\sum_{k = 1}^\infty \frac{f_k(s)}{s + \nu k}
\]
with $f_k(s)$ an entire function of $s$ such that $f_k(-\nu k) =
\nu \Big( \frac{\tau \beta}{\alpha} \Big)^k$. In particular,
$\zeta_{\mathrm{sing}}(s,\Delta_L)$ has possible poles at each $s
= - \nu k$ with the residue equal to
\[
\mathrm{Res}_{s = - \nu k} \zeta_{\mathrm{sing}}(s ,\Delta_L) =
\frac{\sin (\pi (-\nu k))}{\pi} \nu \Big( \frac{\tau
\beta}{\alpha} \Big)^k = - \frac{\nu \sin (\pi \nu k)}{\pi} \Big(
\frac{\tau \beta}{\alpha} \Big)^k,
\]
which is the main result of \cite{FMP} (see equation (7.11) of
loc.\ cit.).

Assume now that $\la = - \frac14$. In this case,
\[
p(x,y) := \det
\begin{pmatrix} \alpha & \beta \\
x & 1 \end{pmatrix} = \alpha - \beta \, x = \alpha\Big( 1 -
\frac{\beta}{\alpha} x \Big),
\]
where we assume that $\alpha,\beta \ne 0$ (the $\alpha = 0$ or
$\beta = 0$ cases can be handled easily), and we write $p(x,y)$ as
in \eqref{pxy}. Forming the power series \eqref{clx}, we see that
\[
\log \Big( 1 - \frac{\beta}{\alpha} x \Big) = \sum_{\ell =
1}^\infty c_{\ell, 0}\, x^\ell y^{2 \cdot 0}\ , \qquad c_{\ell, 0}
= - \frac1{\ell}\Big( \frac{\beta}{\alpha} \Big)^\ell.
\]
Using the definitions \eqref{ellxi} and \eqref{Spl} for $p_\xi$,
$\ell_{\xi}$ (there is only one ``$\xi$" in the present situation,
$\xi = 0$), $\mathscr{\pP}$, and $\mathscr{\lL}$, we immediately
see that $p_{0}$ never exists so $\mathscr{\pP} = \varnothing$,
while
\[
\ell_{0} = \min \{\ell > 0 \, |\, c_{\ell, 0} \ne 0\} = 1 \ ,
\quad \mathscr{\lL} = \{0\}.
\]
Therefore, by Theorem \ref{thm-main},
\[
\zeta_{\mathrm{sing}}(s,\Delta_L) = \frac{\sin (\pi s)}{\pi}
\bigg\{ - e^{-2 s (\log 2 - \gamma)} \log s + g_0(s) \log s
\bigg\},
\]
$g_0(s)$ is an entire function of $s$ such that $g_0(s) =
\mathcal{O}(s)$.  Hence $\zeta(s,\Delta_L)$ has a \emph{genuine}
logarithmic singularity at $s = 0$. This corrects unfortunate
errors from the beautiful paper \cite{FMP} (and \cite{MoE99}),
which states that $\zeta(s,\Delta_L)$ has the ``usual" meromorphic
structure.\footnote{The error in \cite{FMP} occurs in equation
(A13) where certain antiderivatives (specifically, $x N_1(x)$ and
$x^2 N_1(x)^2$) were accidentally set equal to zero at $x = 0$.}
When $\beta = 0$, one can easily see that we still have a
logarithmic singularity at $s = 0$, and when $\alpha = 0$, one can
easily check that there is only the ``regular" part
$\zeta_{\mathrm{reg}}(s,\Delta_L)$ and no ``singular" part; in
fact, the case $\alpha = 0$ corresponds to the Friedrichs
extension (see \cite{BS1}); thus we can see that
$\zeta(s,\Delta_L)$ has a logarithmic singularity for all
extensions except the Friedrichs.
\end{example}

\begin{example} \label{eg-Lap} \em (The Laplacian on $\R^2$)
If $\Delta$ is the Laplacian on a compact region in $\R^2$, then
as we saw before in Section \ref{ssec-egR2}, $A_\Gamma$ has a
$-\frac14$ eigenvalue of multiplicity one and no eigenvalues in
$(- \frac14, \frac34)$. Therefore, the exact same argument we used
in the $\la = - \frac14$ case of the previous example shows that
$\zeta(s,\Delta_L)$ has a logarithmic singularity for all
extensions except the Friedrichs.
\end{example}

\begin{example} \label{eg-countk}
\em Consider now the case of a regular singular operator $\Delta$
over a compact manifold and suppose that $A_\Gamma$ has two
eigenvalues in $[-\frac14,\frac34)$, the eigenvalue $- \frac14$
and another eigenvalue $- \frac14 < \la < \frac34$, both of
multiplicity one. In this case, $V = \C^4$, therefore Lagrangians
$L \subset \C^4$ are determined by $2 \times 2$ matrices $\Aa$ and
$\Bb$. Consider the specific examples
\[
{ \Aa = \begin{pmatrix} 0 & 1 \\
-1 & 0
\end{pmatrix} \ , \quad  \Bb = \Id  .}
\]
Then with $\nu := \sqrt{\la + \frac14}$ and $\tau :=
2^{2\nu}\frac{ \Gamma(1 + \nu)}{ \Gamma(1 - \nu)}$, we have
\[
{ p(x,y) := \det
\begin{pmatrix} 0 & 1 & 1 & 0 \\
-1 & 0 & 0 & 1 \\ x & 0 & 1 & 0 \\ 0 & \tau \, y^{2 \nu} & 0 & 1
\end{pmatrix} = 1 + \tau\, x \, y^{2 \nu}} .
\]
Forming the power series \eqref{clx}, we see that
\[
{ \log \Big( 1 + \tau\, x \, y^{2 \nu} \Big) = \sum_{k = 1}^\infty
\frac{(-1)^{k-1}}{k} \Big( \tau\, x \, y^{2 \nu} \Big)^k = \sum_{k
= 1}^\infty c_{k, \nu k}\, x^k y^{2 \nu k}, }
\]
where ${ c_{k, \nu k} = (-1)^{k-1} \frac{\tau^k}{k}}$. Using the
definitions \eqref{ellxi} and \eqref{Spl} for $p_\xi$,
$\ell_{\xi}$, $\mathscr{\pP}$, and $\mathscr{\lL}$, we immediately
see that $p_{\nu k}$ never exists so $\mathscr{\pP} =
\varnothing$, while
\[
\ell_{\nu k} = \min \{\ell > 0\, |\, c_{\ell , \nu k} \ne 0\} = k
\ , \quad \mathscr{\lL} = \{\nu \, k\, |\, k \in \N\}.
\]
Therefore, by Theorem \ref{thm-main},
\[
\zeta_{\mathrm{sing}}(s,\Delta_L) = \frac{\sin (\pi s)}{\pi}
\bigg\{ - e^{-2 s (\log 2 - \gamma)} \log s + \sum_{k = 1}^\infty
g_k(s) \log (s + \nu k) \bigg\},
\]
with $g_k(s)$ an entire function of $s$ such that near $s = - \nu
k$,
\[
{ g_k(s) = (-1)^k \frac{\tau^k 2^{k} \nu}{(k - 1)!} (s + \nu k)^{k
- 1} + \mathcal{O}((s + \nu k)^{k}) . }
\]
In particular, $\zeta_{\mathrm{sing}}(s,\Delta_L)$ has
\emph{countably} many logarithmic singularities!
\end{example}

\begin{example} \em With the same situation as considered in Example
\ref{eg-countk}, consider
\[
{ \Aa = \begin{pmatrix} - 1 & 1 \\
0 & 0 \end{pmatrix}} \ , \quad \Bb = \begin{pmatrix} 0 & 0 \\
1 & -1
\end{pmatrix} .
\]
Then with $\nu := \sqrt{\la + \frac14}$ and $\tau :=
2^{2\nu}\frac{ \Gamma(1 + \nu)}{ \Gamma(1 - \nu)}$, we have
\[
{ p(x,y) := \det
\begin{pmatrix} -1 & 1 & 0 & 0 \\
0 & 0 & 1 & -1 \\ x & 0 & 1 & 0 \\ 0 & \tau \, y^{2 \nu} & 0 & 1
\end{pmatrix} = x - \tau\, y^{2 \nu} = x \Big(1 - \tau \, x^{-1} y^{2 \nu}
\Big).}
\]
Forming the power series \eqref{clx}, we see that
\[
{ \log \Big( 1 - \tau \, x^{-1} \, y^{2 \nu} \Big) = \sum_{k =
1}^\infty \frac{(-1)^{k-1}}{k} \Big(-\tau\, x^{-1} y^{2 \nu}
\Big)^k = \sum_{k = 1}^\infty c_{-k, \nu k}\, x^{-k} y^{2 \nu k},}
\]
where ${ c_{-k, \nu k} = - \frac{\tau^k}{k}}$. Using the
definitions \eqref{ellxi} and \eqref{Spl} for $p_\xi$,
$\ell_{\xi}$, $\mathscr{\pP}$, and $\mathscr{\lL}$, we immediately
see that $\ell_{\nu k}$ never exists so $\mathscr{\lL} =
\varnothing$, while
\[
p_{\nu k} = \min \{\ell \leq 0\, |\, c_{-\ell , \nu k} \ne 0\} =
-k \ , \quad \mathscr{\pP} = \{\nu \, k\, |\, k \in \N\} .
\]
Therefore, by Theorem \ref{thm-main},
\[
\zeta_{\mathrm{sing}}(s,\Delta_L) = \frac{\sin (\pi s)}{\pi}
\sum_{k = 1}^\infty \frac{f_k(s)}{(s + \nu k)^{k + 1}}
\]
with $f_k(s)$ an entire function of $s$ such that $f_k(-\nu k) =
(-1)^k \frac{\tau^k k! \nu }{2^{k}}$. In particular,
$\zeta_{\mathrm{sing}}(s,\Delta_L)$ has poles of
\emph{arbitrarily} large order!
\end{example}

\begin{example} \label{eg-count3k}
\em Consider one last example, the case of a regular singular
operator $\Delta$ over a compact manifold such that $A_\Gamma$ has
three eigenvalues in $[-\frac14,\frac34)$, the eigenvalue $-
\frac14$ with multiplicity two and another eigenvalue $- \frac14 <
\la < \frac34$ of multiplicity one. In this case, $V = \C^6$ and
Lagrangians $L \subset \C^6$ are determined by $3 \times 3$
matrices $\Aa$ and $\Bb$. Consider the specific examples
\[
{ \Aa = \begin{pmatrix} 0 & 1 & -1 \\
1 & 0 & 0 \\
1 & 0 & 0
\end{pmatrix}} \ , \quad \Bb = \Id.
\]
Then with $\nu := \sqrt{\la + \frac14}$ and $\tau :=
2^{2\nu}\frac{ \Gamma(1 + \nu)}{ \Gamma(1 - \nu)}$, we have {
\begin{align*} p(x,y) & := \det
\begin{pmatrix} 0 & 1 & -1 & 1 & 0 & 0 \\
1 & 0 & 0 & 0 & 1 & 0 \\ 1 & 0 & 0 & 0 & 0 & 1\\
x & 0 & 0 & 1 & 0 & 0 \\ 0 & x & 0 & 0 & 1 & 0\\
0 & 0 & \tau \, y^{2 \nu} & 0 & 0 & 1
\end{pmatrix}\\ & \qquad = - x + \tau\, y^{2 \nu} - \tau\, x^2 \, y^{2 \nu} =
- x \Big( 1 + \big( x - x^{-1} \big) \tau\, y^{2 \nu} \Big).
\end{align*}}
Forming the power series \eqref{clx}, we see that
\begin{multline} \label{ex2logs} \log \Big( 1 + \big( x - x^{-1}
\big) \tau\, y^{2 \nu} \Big) = \sum_{k = 1}^\infty
\frac{(-1)^{k-1}}{k} \tau^k \big( x - x^{-1}
\big)^k y^{2 \nu k} \\
= \sum_{k = 1}^\infty \sum_{j = 0}^k \frac{(-1)^{j-1}}{k} \tau^k
\binom{k}{j} \, x^{2 j - k} y^{2 \nu k} = \sum_{\ell, k}\,
c_{\ell, \nu k} \, x^\ell y^{2 \nu k},
\end{multline}
where for each $k$, $\ell$ runs from $-k$ to $k$. Using the
definitions \eqref{ellxi} and \eqref{Spl} for $p_\xi$,
$\ell_{\xi}$, $\mathscr{\pP}$, and $\mathscr{\lL}$, we immediately
see that
\[
p_{\nu k} = \min \{\ell \leq 0 \, |\, c_{\ell , \nu k} \ne 0\} =
-k,
\]
\[
\ell_{\nu k} = \min \{\ell > 0\, |\, c_{\ell , \nu k} \ne 0\} =
\begin{cases} 1 & \text{if $k$ is odd,}\\ 2 & \text{if $k$ is even,}
\end{cases},
\]
and $\mathscr{\pP} = \mathscr{\lL} = \{\nu \, k\, |\, k \in \N\}$.
Therefore, by Theorem \ref{thm-main},
\begin{multline}
\zeta_{\mathrm{sing}}(s,\Delta_L) = \frac{\sin (\pi s)}{\pi}
\bigg\{ - e^{-2 s (\log 2 - \gamma)} \log s + \sum_{k = 1}^\infty
\frac{f_k(s)}{(s + \nu k )^{k + 1}}\\
+ \sum_{k = 1}^\infty g_k(s) \log (s + \nu k) \bigg\},
\end{multline}
where $f_k(s)$ and $g_k(s)$ are entire functions of $s$ such that
$f_k(-\nu k) = (-1)^k \frac{ \tau^k k! \nu }{2^{k}}$ and
\[
{ g_k(s) = 2 \nu (-1)^{m+1} \tau^k \binom{k}{m+1} \times
\begin{cases} 1 + \mathcal{O}((s + \nu k))\ \ \text{if $k = 2 m +
1$ is odd,}\\ 2 (s + \nu k) + \mathcal{O}((s + \nu k)^{2})\ \
\text{if $k = 2 m$ is even.}
\end{cases}}
\]
In particular, $\zeta_{\mathrm{sing}}(s,\Delta_L)$ has poles of
\emph{arbitrarily} high orders and in addition to a logarithmic
singularity at the origin, \emph{countably} many logarithmic
singularities at the same locations of the poles!
\end{example}

\subsection{Special Lagrangian subspaces} \label{ssec-special}
Theorem \ref{thm-main} simplifies considerably when the conditions
given by the Lagrangian $L$ over $- \frac14$ eigenspaces are
separated from the conditions over $\la$ eigenspaces with $\la$ in
$(-\frac14,\frac34)$. We shall call a Lagrangian subspace $L
\subset V$ \emph{decomposable} if $L = L_0 \oplus L_1$ where $L_0$
is an arbitrary Lagrangian subspace of $\bigoplus_{\la_\ell =
-\frac14} E_\ell\oplus E_\ell$ and $L_1$ is an arbitrary
Lagrangian subspace of $\bigoplus_{-\frac14<\la_\ell<\frac34}
E_\ell\oplus E_\ell$. As described in Proposition
\ref{p:relation}, the Lagrangian subspace $L_0$ is determined by
two $q_0 \times q_0$ matrices $\Aa_0$, $\Bb_0$ where $q_0 = \dim
L_0$, that is, the multiplicity of the eigenvalues $\la_\ell =
-\frac14$. Similarly, the Lagrangian subspace $L_1$ is determined
by two $q_1 \times q_1$ matrices $\Aa_1$, $\Bb_1$ where $q_1 =
\dim L_1$, that is, the multiplicity of the eigenvalues $\la_\ell$
with $-\frac14<\la_\ell<\frac34$.

Two polynomials which are explicitly determined by the matrices
$\Aa_0,\Bb_0$ and $\Aa_1,\Bb_1$ play peculiar r\^{o}les in the
statement of our result. First, consider the polynomial $p_0(z)$
in the single variable $z$ defined by
\begin{equation} \label{p0(z)}
p_0(z) := \det \begin{pmatrix} \Aa_0 & \Bb_0 \\
\Id_{q_0} & (\log 2 - \gamma - z) \Id_{q_0} \end{pmatrix} .
\end{equation}
Using the definition of determinant, it is easy to see that
$p_0(z)$ is a polynomial of degree at most $q_0$ in $z$. Since the
degree of $p_0'(z)$ is one less than the degree of $p_0(z)$, we
can write
\begin{equation} \label{p0'p0}
\frac{p_0'(z)}{p_0(z)} = \sum_{k = 1}^\infty \frac{\beta_k}{z^k} \
, \qquad \beta_k \in \C,
\end{equation}
where the series on the right is absolutely convergent for $|z|$
sufficiently large. Second, consider the polynomial in
\eqref{polyp} (and \eqref{pxy}) using $\Aa_1$ and $\Bb_1$ in place
of $\Aa$ and $\Bb$:
\begin{equation}\label{p1}
p_1(y) := \det
\begin{pmatrix} \Aa_1 & \Bb_1\\
\begin{array}{ccc}
\hspace{-.5em} \tau_1\, y^{2 \nu_1} & 0 & 0\\
0 & \ddots & 0\\
0 & 0 & \tau_{q_1}\, y^{2 \nu_{q_1}}
\end{array}
&  \Id_{q_1} \hspace{-.3em}
\end{pmatrix} = a_{\alpha_0}\, y^{2 \alpha_0} \Big( 1 +
\sum b_\beta y^{2 \beta} \Big),
\end{equation}
where the $\beta$'s are positive. Then as in \eqref{clx}, write
\begin{equation} \label{cxi}
\log \Big( 1 + \sum  b_{\beta} y^{2 \beta} \Big) = \sum c_{\xi} \,
y^{2 \xi} ,
\end{equation}
and let $\mathscr{\pP} := \{ \xi \, |\, c_\xi \ne 0\}$. Then
Theorem \ref{thm-main} simplifies to

\begin{theorem} \label{main thm} For an arbitrary decomposable
Lagrangian $L \subset V$, the $\zeta$-function $\zeta(s,\Delta_L)$
has the following form:
\begin{equation}\label{zeta-dec}
\z(s,\Delta_L) = \z_{\mathrm{reg}}(s,\Delta_L) +
\zeta_{\mathrm{sing}}(s,\Delta_L) ,
\end{equation}
where $\zeta_{\mathrm{reg}}(s,\Delta_L)$ has the ``regular" poles
at the ``usual" locations $s = \frac{n - k}{2}\notin -\N_0$ for $k
\in \N_0$ and at $s = 0$ if $\dim \Gamma > 0$, and where
$\zeta_{\mathrm{sing}}(s,\Delta_L)$ has the following expansion:
\begin{equation}\label{e:fs}
\zeta_{\mathrm{sing}}(s,\Delta_L) =  - \frac{ \sin (\pi
s)}{\pi}f(s) \log s \ +\ \frac{\sin (\pi s)}{\pi} \! \! \sum_{\xi
\in \mathscr{P}} \frac{f_\xi(s)}{s + \xi},
\end{equation}
where $f(s)$ is the entire function defined explicitly by $f(s) =
\sum_{ k = 1 }^\infty \beta_k \frac{(- 2 s)^{k-1}}{(k-1)!}$, and
the $f_\xi(s)$'s are entire functions of $s$ such that
\[
f_\xi(-\xi) = - {c_{\xi} \,  \xi}.
\]
\end{theorem}

For certain types of Lagrangians, the formula for $f(s)$ becomes
very simple. We shall call the Lagrangian $L_0$ \emph{split-type}
if it can be written as $\bigoplus_{\la_\ell = - \frac14}
L_{\ell}$ with $L_\ell$ a Lagrangian subspace of $E_\ell \oplus
E_\ell$. As explained in Proposition \ref{p:Ltheta}, each
component $L_\ell$ of $L_0$ is determined by an angle
$\theta_\ell\in [0,\pi)$. Moreover, in this case, the coefficients
$\beta_k$ in \eqref{p0'p0} are given by (as follows from Corollary
\ref{cor-m=1caseint})
\begin{equation} \label{betak}
\beta_k = \sum_{\theta_\ell \ne \frac{\pi}{2}}
\kappa_\ell^{k-1},\qquad k = 1,2,3,4,\ldots. \quad \text{(for
split-type $L_0$)}
\end{equation}
where $\kappa_\ell = \log 2 - \gamma - \tan \theta_\ell$ with
$\theta_\ell$ the angle defining $L_\ell$ in $L_0$ and $\gamma$ is
the Euler-Mascheroni constant. Then when $L_0$ is of split-type,
we have

\begin{theorem} \label{main thm2} For a decomposable
Lagrangian $L \subset V$ such that $L_0$ is of split-type,
$\zeta(s,\Delta_L)$ has the form as in \eqref{zeta-dec} with $f(s)
= \sum_{\theta_\ell \ne \frac{\pi}{2}} e^{-2 s \kappa_\ell}$ in
\eqref{e:fs}.
\end{theorem}

\begin{example} \label{eg-LapR2} \em
Consider the case when $A_\Gamma$ has exactly one eigenvalue in
$[-\frac14,\frac34)$, the eigenvalue $-\frac14$. Then Theorem
\ref{main thm2} shows the result stated in (\ref{add1}).

\end{example}

\subsection{Unusual resolvent and heat kernel expansions}

Besides establishing exotic $\zeta$-expansions, we also derive
equally exotic resolvent and heat kernel expansions.

\begin{theorem} \label{thm-restracegen}
Let $\Lambda \subset \C$ be any sector (solid angle) not
intersecting the positive real axis and choose $N \in \N$ with $N
\geq \frac{n}{2}$. Then for an \textbf{arbitrary} Lagrangian $L$,
as $|\la| \to \infty$ with $\la \in \Lambda$ we have
\begin{multline} \label{resexpgen}
\Tr (\Delta_L - \lambda)^{-N-1} \,  \sim \, \sum_{k=0}^\infty a_k
\, (- \lambda)^{\frac{n-k}{2} - N-1}  + b \,  (- \lambda)^{- N -
1} \, \log (- \lambda)\\
\quad \quad +\frac{1}{N!} \frac{d^N}{d \la^N} \bigg\{ \frac{q_0 -
j_0}{(-\la) (\log (-\la) - 2 \widetilde{\gamma})} \bigg\} \\ -
\frac{1}{N!}\frac{d^{N+1}}{d \la^{N+1}} \bigg\{ \sum 2^\ell
c_{\ell \xi}\, (-\la)^{- \xi} \Big(2 \widetilde{\gamma} - \log
(-\la) \Big)^{-\ell} \bigg\}
\end{multline}
where the $a_k$ and $b$ coefficients are independent of $L$, the
$c_{\ell \xi}$'s are the coefficients in \eqref{clx}, and
$\widetilde{\gamma} = \log 2 -\gamma$.
\end{theorem}

From the explicit formula \eqref{resexpgen} and from the binomial
theorem for $\ell>0$:
\begin{align} \label{resexample}
\Big(2 \widetilde{\gamma} - \log (-\la) \Big)^{-\ell} & = (- \log
(-\la))^{-\ell} \bigg(1 - \frac{2 \widetilde{\gamma}}{\log (-\la)}
\bigg)^{-\ell} \\ \notag & = (- \log (-\la))^{-\ell} \sum_{j =
0}^\infty \binom{-\ell}{j} \frac{ (-2 \widetilde{\gamma})^j}{
(\log (-\la))^j},
\end{align}
it is obvious that when $A_\Gamma$ has $-\frac14$ eigenvalues, the
resolvent trace expansion has, in general, $\log(-\la)$ terms of
\emph{arbitrarily} high multiplicity and inverse powers $(\log
(-\la))^{-1}$ with \emph{infinite} multiplicity! This phenomenon
is new and even for pseudodifferential operators on compact
manifolds, with or without boundary, and even conic, ``regular"
(not inverse powers of) $\log (- \la)$ terms occur with at most
multiplicity two \cite{G,GL,GL2,GrG-HaL02,GrG-SeR95,LoP01,LoPR02}.
See \cite{GKM,K,SS} for studies of resolvents for closed
extensions of general cone operators in the sense of Schulze
\cite{Sz91}. Here is a concrete example illustrating this
discussion:

\begin{example}
\em For the self-adjoint extension $\Delta_L$ considered in
Example \ref{eg-count3k}, from the explicit formula
\eqref{ex2logs}, we immediately get
\begin{multline*}
\Tr (\Delta_L - \lambda)^{-N-1} \,  \sim \, \sum_{k=0}^\infty a_k
\, (- \lambda)^{\frac{n-k}{2}  - N-1}  + b \,  (- \lambda)^{- N -
1} \, \log (- \lambda)\\
+\frac{1}{N!} \frac{d^N}{d \la^N} \bigg\{ \frac{1}{(-\la) (\log
(-\la) - 2 \widetilde{\gamma})} \bigg\} \\ -\frac{1}{N!}
\frac{d^{N+1}}{d \la^{N+1}} \bigg\{ \sum_{k = 1}^\infty \sum_{j =
0}^k \frac{2^{2j - k} (-1)^{ j - 1}}{k} \tau^k \binom{k}{j} \,
(-\la)^{- 2 \nu k} \big(2 \widetilde{\gamma} - \log (- \la) \big)
^{k - 2 j} \bigg\}.
\end{multline*}
In this very simple example, we see unusual powers $(-\la)^{- 2\nu
k - N - 1}$ (after taking $N+1$ derivatives) and log terms
$\log(-\la)$ of \emph{arbitrarily} high multiplicity (each unusual
power $(-\la)^{- 2\nu k - N - 1}$ with a log term of highest power
$(\log (-\la))^k$), and inverse powers $(\log (-\la))^{-1}$ with
\emph{infinite} multiplicity because of the formula
\eqref{resexample}.
\end{example}

When $L$ is decomposable, the last two terms in \eqref{resexpgen}
can be made very explicit.

\begin{theorem} \label{thm-restrace}
Let $\Lambda \subset \C$ be any sector (solid angle) not
intersecting the positive real axis and choose $N \in \N$ with $N
\geq \frac{n}{2}$. Then for an arbitrary decomposable Lagrangian
$L$, as $|\la| \to \infty$ with $\la \in \Lambda$ we have
\begin{align}
\Tr (\Delta_L - \lambda)^{-N-1} & \sim \sum_{k=0}^\infty a_k \, (-
\lambda)^{\frac{n-k}{2}  - N-1}  + b \,  (- \lambda)^{- N - 1} \,
\log (- \lambda) \label{resexp} \\   - \frac{1}{N!}&
\frac{d^{N+1}}{d \la^{N+1}} \bigg\{ \sum_{\xi \in \mathscr{\pP}}
c_{\xi}\, (-\la)^{- \xi} \bigg\} +\frac{1}{N!} \frac{d^N}{d \la^N}
\left\{ (-\la)^{-1} \sum_{k = 1}^\infty \frac{2^{k-1}
\beta_k}{(\log (-\la))^{k}} \right\}  , \notag
\end{align}
where the $a_k$ and $b$ coefficients are independent of $L$, the
$\beta_k$'s are the coefficients in \eqref{p0'p0}, and the
$c_{\xi}$'s are the coefficients in \eqref{cxi}.
\end{theorem}

In the case when $L_0$ is of split-type, the second-to-last term
in \eqref{resexp} can be made even more explicit because of the
formula \eqref{betak} for $\beta_k$. We also prove a corresponding
heat kernel expansion.

\begin{theorem} \label{t:traceM}
For an \textbf{arbitrary} Lagrangian $L$, the heat kernel $e^{-t
\Delta_L}$ has the following trace expansion as $t \to 0$:
\begin{align*}
\Tr ( e^{-t \Delta_{L}})\, \sim \,  & \sum_{k=0}^\infty
\widetilde{a}_k\, t^{\frac{-n+k}{2}} + b\log t
+\sum_{k=0}^\infty \widetilde{b}_k (\log t)^{-1-k} \\& +
\sum_{\xi\in \mathscr{\pP}}\sum^{|p_\xi|+1}_{k=0}
\widetilde{c}_{\xi k}\, t^{\xi} (\log t)^{k} + \sum_{\xi\in
\mathscr{\lL}}\sum_{k=0}^\infty \widetilde{d}_{\xi k}\, t^{\xi}\,
(\log t)^{-\ell_\xi-k} ,
\end{align*}
with $\widetilde c _{10} = 0$ and $\widetilde c_{\xi (|p_\xi| +1)}
=0$ for $\xi \notin \N_0$.
\end{theorem}
Thus, the heat trace expansion, in general, has powers of $\log t$
with finite multiplicity and inverse powers $(\log t)^{-1}$ with
\emph{infinite} multiplicity. The $\widetilde{c}_{\xi k}$ and
$\widetilde{d}_{\xi k}$ coefficients can be expressed in terms of
the coefficients in the resolvent expansion \eqref{resexpgen} but
not so explicitly. For decomposable Lagrangians we have

\begin{theorem} \label{t:traceMde}
For an arbitrary decomposable Lagrangian $L$, the heat kernel
$e^{-t \Delta_L}$ has the following trace expansion as $t \to 0$:
\[
\Tr ( e^{-t \Delta_{L}})\, \sim \, \sum_{k=0}^\infty
\widetilde{a}_{k}\, t^{\frac{-n+k}{2}}+\widetilde{b}\, \log t +
\sum_{\xi \in \mathscr{\pP}} \widetilde{c}_\xi \, t^{\xi} +
\sum_{k = 1}^\infty \widetilde{d}_k(\log t)^{-k}.
\]
\end{theorem}

\section{Hermitian symplectic formalism of self-adjoint extensions}
\label{s:SAE}

To orient the reader to the various terminologies used throughout
this paper, in this section we briefly review the classical theory
of the Hermitian forms characterization of self-adjoint
extensions. For more on this viewpoint, see \cite{HarM00,HarMI00}
and see \cite{BExP-SeP89} for applications of self-adjoint
extensions to quantum physics.

\subsection{The maximal domain}
An elegant observation due to Gelfand around $1971$ (according to
Novikov \cite[p.\ 1]{INovS99}) to analyze self-adjoint extensions
of $\Delta$ is to find ``maximal" domains $\mathfrak{D} \subset
\dom_\max(\Delta )$ that make the Hermitian quadratic form
\[
\langle \Delta \phi ,\psi \rangle - \langle \phi , \Delta \psi
\rangle,\qquad \phi ,\psi \in \mathfrak{D},
\]
vanish. The use of Hermitian forms as a tool to analyze
self-adjoint extensions of (one-dimensional) regular singular
operators goes back at least to Kochube{\u\i} in the late $70$'s
\cite{KocA79} (cf.\ \cite{KocA90,KocA91}) and Pavlov \cite{PavB87}
in the late $80$'s.  Mooers \cite{MoE99} (cf.\ \cite{Ch1,Ch2,LMP},
and Gil and Mendoza \cite{GM} for general cone operators) used
this technique to give a complete description of the self-adjoint
extensions of the regular singular operator $\Delta$ in terms of
Lagrangian subspaces in the eigenspaces of $A_\Gamma$ with
eigenvalues in the interval $[- \frac14, \frac34)$. To describe
these extensions, we first describe $\mathrm{dom}_{\max}(\Delta)$.
Recall that $ \Delta|_{\mathcal{U}}=- \partial^2_r + \frac{1}{r^2}
A_\Gamma$ (see \eqref{D}), where $\mathcal{U}\cong
[0,\varepsilon)_r\times \Gamma$ denotes a tubular neighborhood of
$\Gamma$ and $A_\Gamma$ is a Laplace-type operator over $\Gamma$
such that $A_\Gamma \geq -\frac14$. Recall that
\[
-\frac14 = \underbrace{\lambda_1 = \lambda_2 = \cdots =
\lambda_{q_0}}_{ = - \frac14} < \underbrace{\lambda_{q_0+1} \leq
\lambda_{q_0 + 2} \leq \cdots \leq \la_{q_0 + q_1}}_{- \frac14 <
\la_\ell < \frac34}
\]
denotes the eigenvalues of $A_\Gamma$ in $[- \frac14, \frac34)$
with corresponding orthonormal eigenvectors $\{\phi_\ell\}$. Then,
as shown by Cheeger \cite{Ch1,Ch2} (cf.\ Mooers \cite{MoE99}), we
have

\begin{proposition} \label{p:maxdom}
A section $\phi \in L^2(M,E)$ is in
\[
\dom_{\max}(\Delta) := \left\{ \phi \in L^2(M,E)\ |\ \Delta \phi
\in L^2(M,E) \right\},
\]
where ``$\Delta \phi \in L^2(M,E)$" is in the distributional
sense, if and only if $\phi$ is in $H^2$ away from the boundary
$\Gamma$, and near $\Gamma$ we can write
\begin{align} \label{phiasymp1}
\phi =& \sum_{\ell = 1}^{q_0} \left\{ {c_\ell^+(\phi)} \, r^{
\frac12} \phi_\ell + {c_{\ell}^-(\phi)} \,
r^{ \frac12}\log r \phi_\ell \right\}\\
 &\qquad + \sum_{\ell = 1}^{q_1} \left\{ \frac{c_{q_0 + \ell}^+(\phi)}{\sqrt{2
\nu_\ell}} \, r^{\nu_\ell + \frac12} \phi_{q_0 + \ell} +
\frac{c_{q_0 + \ell}^-(\phi)}{\sqrt{2 \nu_\ell}} \, r^{-\nu_\ell +
\frac12} \phi_{q_0 + \ell} \right\} + \widetilde{\phi} ,\notag
\end{align}
where the $c_\ell^\pm(\phi)$'s are constants, $\nu_\ell :=
\sqrt{\lambda_{q_0 + \ell} + \frac14} > 0$, $\widetilde{\phi} \in
H^2$ and $\widetilde{\phi}=\mathcal{O}(r^{\frac32})$.
\end{proposition}

\subsection{Self-adjoint extensions}

Given a domain $\mathfrak{D} \subset \mathrm{dom}_{\max}(\Delta)$,
we say that
\[
\Delta_{\mathfrak{D}} := \Delta : \mathfrak{D} \to L^2(M,E)
\]
is \emph{self-adjoint} if
\[
\{ \psi \in \mathrm{dom}_\max(\Delta)\ |\ \langle \Delta \phi ,
\psi \rangle = \langle \phi , \Delta \psi \rangle \ \ \text{for
all}\ \ \phi \in \mathfrak{D} \} = \mathfrak{D}.
\]
Simply put: $\Delta$ is ``maximally" symmetric on $\mathfrak{D}$
in the sense that $\Delta$ is symmetric on $\mathfrak{D}$ and
adding any elements to $\mathfrak{D}$ will destroy this symmetry.
In this sense, $\mathfrak{D}$ is a ``maximal" domain that makes
the Hermitian form $\langle \Delta \phi ,\psi \rangle - \langle
\phi , \Delta \psi \rangle$ vanish.

Now, an integration by parts argument (see \cite[Prop.\
2.4]{MoE99}) shows that
\begin{multline} \label{sympform}
\langle \Delta \phi, \psi\rangle - \langle  \phi, \Delta \psi
\rangle\ \ = -  \sum_{\ell = 1}^{q_0} \! \! \left( c_\ell^+(\phi)
\overline{c_\ell^-(\psi)} - c_\ell^-(\phi)
\overline{c_\ell^+(\psi)} \right) \\ + \! \! \! \! \sum_{\ell =
q_0 + 1}^q \! \! \left( c_\ell^+(\phi) \overline{c_\ell^-(\psi)} -
c_\ell^-(\phi) \overline{c_\ell^+(\psi)} \right) ,\qquad \text{for
all}\ \phi, \psi \in \dom_{\max}(\Delta).
\end{multline}
To put this in a symplectic framework, we define $\psi^+_\ell : =
r^{\frac12 + \nu_\ell}\phi_\ell$ for
$0 \leq \nu_\ell <1$, and
\[
\psi^-_\ell := \begin{cases} r^{\frac12 - \nu_\ell} \phi_\ell &
\text{for $0 < \nu_\ell <1$,}\\
r^{\frac12} \log r\, \phi_\ell & \text{for $\nu_\ell = 0$} ,
\end{cases}
\]
and define
\[
V : = \bigoplus_{-\frac14 \leq \la_\ell < \frac 34} E_\ell^+\oplus
E^-_\ell
\]
where $ E_\ell^{\pm} := \langle\, \psi_\ell^\pm\, \rangle$ with
$\langle\, \psi_\ell^\pm\, \rangle :=
\mathrm{span}_\C\{\psi_\ell^\pm\}$. We endow $V$ with the
symplectic structure $\omega: V \times V \to \C$ defined by
\begin{equation} \label{symplectic2.5}
\omega(\psi^\pm_\ell,\psi^\mp_\ell) = \begin{cases} { \mp 1} &
\hbox{
when } \nu_\ell = 0 \\
{ \pm 1} &  \hbox{ when } \nu_\ell \ne 0; \end{cases}\qquad
\omega(\psi^\pm_\ell,\psi^\pm_j)  = 0 \ \ \hbox{otherwise},
\end{equation}
and extending to $V \times V$ linearly in the first factor and
conjugate linearly in the second factor. Then according to
\eqref{phiasymp1} and \eqref{sympform}, we have (cf.\ \cite[Prop.\
2.4]{MoE99}):

\begin{lemma} \label{lem-omega}
$\omega : V \times V \to \C$ is a Hermitian symplectic form and
for any $\phi,\psi \in \mathrm{dom}_{\max}(\Delta)$, we have
\[
{ \langle \Delta \phi, \psi\rangle - \langle  \phi, \Delta \psi
\rangle } = \omega( \vec{\phi}, \vec{\psi})
\]
where $\vec{\phi} := \sum \big\{ c^+_\ell (\phi) \psi^+_\ell +
c^-_\ell (\phi) \psi^-_\ell \big\} \in V$ and $\vec{\psi}$ is
defined similarly.
\end{lemma}

A subspace $L \subset V$ is called \emph{Lagrangian} if
\[
\{ w \in V\ |\ \omega( v , w) = 0 \ \ \text{for all}\ \ v \in L \}
= L.
\]
We can now prove:

\begin{theorem} \label{thm-sae}
Self-adjoint extensions of $\Delta$ are in one-to-one
correspondence with Lagrangian subspaces of $V$ in the sense that
given any Lagrangian subspace $L \subset V$, defining
\[
\mathrm{dom}_L(\Delta) := \{ \phi \in \mathrm{dom}_\max(\Delta) \
|\ \vec{\phi} \in L\},
\]
the operator
\[
\Delta_L := \Delta : \mathrm{dom}_L(\Delta) \to L^2(M,E)
\]
is self-adjoint and any self-adjoint extension of $\Delta$ is of
the form $\Delta_L$ for some Lagrangian $L \subset V$.
\end{theorem}
\begin{proof}
By Lemma \ref{lem-omega}, we can write this as:
$\Delta_{\mathfrak{D}}$ is self-adjoint if and only if
\begin{equation} \label{Lagcondi}
\omega(\vec{\phi}, \vec{\psi}) = 0 \ \ \text{for all}\ \ \phi \in
\mathfrak{D} \quad \Longleftrightarrow \quad \psi \in
\mathfrak{D}.
\end{equation}

Suppose that $\Delta_{\mathfrak{D}}$ is self-adjoint and define $L
:= \{ \vec{\phi} \in V\ |\ \phi \in \mathfrak{D}\}$; we shall
prove that $L$ is Lagrangian. Let $w \in L$ and choose $\psi \in
\mathfrak{D}$ such that $\vec{\psi} = w$. Then by
\eqref{Lagcondi}, $\omega(\vec{\phi} , w) = 0$ for all $\phi \in
\mathfrak{D}$. Therefore, $\omega(v,w) = 0$ for all $v \in L$.
Conversely, let $w \in V$ and assume that $\omega(v,w) = 0$ for
all $v \in L$. Choose $\psi \in \mathrm{dom}_\max(\Delta)$ such
that $\vec{\psi} = w$; this can always be done, for if $w = \sum
\big\{ c^+_\ell \psi^+_\ell + c^-_\ell  \psi^-_\ell \big\}$, then
\begin{equation} \label{psiconstruct}
\psi := \rho \left(\sum \big\{ c^+_\ell \psi^+_\ell + c^-_\ell
\psi^-_\ell \big\} \right)
\end{equation}
will do, where $\rho \in C^\infty(M)$ is supported in the tubular
neighborhood $\mathcal{U}$ and equals $1$ near $r = 0$. Then
$\omega(v,w) = 0$ for all $v \in L$ implies that
$\omega(\vec{\phi} , \vec{\psi}) = 0$ for all $\phi \in
\mathfrak{D}$, which by \eqref{Lagcondi}, implies that $\psi \in
\mathfrak{D}$, which further implies that $w = \vec{\psi} \in L$.

Now let $L \subset V$ be Lagrangian; we shall prove that
$\Delta_L$ is self-adjoint, that is, \eqref{Lagcondi} holds. Let
$\psi \in \mathrm{dom}_L(\Delta)$. Then, since $L$ is Lagrangian,
we automatically have $\omega(\vec{\phi} , \vec{\psi}) = 0$ for
all $\phi \in \mathrm{dom}_L(\Delta)$. Conversely, let $\psi \in
\mathrm{dom}_\max(\Delta)$ and assume that $\omega(\vec{\phi} ,
\vec{\psi}) = 0$ for all $\phi \in \mathrm{dom}_L(\Delta)$. By the
construction \eqref{psiconstruct}, given any $v \in L$ we can find
a $\phi \in \mathrm{dom}_\max(\Delta)$ such that $\vec{\phi} = v$.
Therefore, $\omega(\vec{\phi} , \vec{\psi}) = 0$ for all $\phi \in
\mathrm{dom}_L(\Delta)$ implies that $\omega(v , \vec{\psi}) = 0$
for all $v \in L$, which by the Lagrangian condition on $L$,
implies that $\vec{\psi} \in L$. This shows that $\psi \in
\mathrm{dom}_L(\Delta)$ and our proof is complete.
\end{proof}

\subsection{Characterizations of Lagrangian subspaces} \label{sec-decomp}

Kostrykin and Schrader \cite[Lem.\ 2.2]{KoV-ScR99} characterize
all Lagrangian subspaces in complex Euclidean space.

\begin{proposition}\label{p:relation}
A subset $L \subset \C^{2k}$ is Lagrangian (with respect to the
standard Euclidean symplectic form) if and only if it can be
described by a system of equations
\[
L = \left\{ \vec{c} \in \C^{2 k} \, |\, \begin{pmatrix} \Aa & \Bb
\end{pmatrix} \vec{c} = 0 \right\} \subset \C^{2 k},
\]
where $\Aa$ and $\Bb$ are $k \times k$ matrices such that the rank
of $\begin{pmatrix} \Aa & \Bb
\end{pmatrix}$ is $k$ and $\Aa \Bb^* = \Bb \Aa^*$.
\end{proposition}

As seen in \eqref{symplectic2.5}, $V$ can be identified with
$\C^{2q} = \C^{2q_0} \times \C^{2q_1}$, where $q = q_0 + q_1$,
with minus the standard symplectic form on the $\C^{2q_0}$ factor
and the standard symplectic form on the $\C^{2q_1}$ factor. Using
this fact, we prove

\begin{corollary} \label{cor-L}
A Lagrangian subspace $L \subset V$ can be characterized by $q
\times q$ matrices $\Aa$ and $\Bb$ via
\[
L \cong \left\{ \vec{c} \in \C^{2 q} \, |\, \begin{pmatrix} \Aa &
\Bb
\end{pmatrix} \vec{c} = 0 \right\} \subset \C^{2 q},
\]
where $\begin{pmatrix} \Aa & \Bb \end{pmatrix}$ has rank $q$ and
$\Aa' \Bb^*$ is self-adjoint where $\Aa'$ is the matrix $\Aa$ with
the first $q_0$ columns of $\Aa$ multiplied by $-1$.
\end{corollary}
\begin{proof}
Define $T : \C^{2 q} \to \C^{2q}$ by
\[
T(v_1,\ldots, v_{q_0},v_{q_0 + 1} ,\ldots , v_{2q}) =
(-v_1,\ldots,- v_{q_0},v_{q_0 + 1} ,\ldots , v_{2q}).
\]
Then using the definition of $\omega$ in \eqref{symplectic2.5}, it
follows that
\[
\omega(v,w) = \omega_E(T v, Tw) \qquad \text{for all\ }v,w \in
\C^{2q},
\]
where we identify $V$ with $\C^{2q}$ and $\omega_E$ is the
standard Euclidean symplectic form on $\C^{2q}$. Therefore,
\begin{align*}
\big\{ \vec{c} \in \C^{2 q} \, & |\, \begin{pmatrix} \Aa & \Bb
\end{pmatrix} \vec{c} = 0 \big\} \ \ \text{is Lag.\ in\ } (\C^{2q},
\omega) \\
 & \Longleftrightarrow \quad  \big\{ T \vec{c} \in \C^{2 q} \, |\,
\begin{pmatrix} \Aa & \Bb
\end{pmatrix} \vec{c} = 0 \big\} \ \ \text{is Lag.\ in\ } (\C^{2q},
\omega_E) \\
& \Longleftrightarrow \quad  \big\{ \vec{c} \in \C^{2 q} \, |\,
\begin{pmatrix} \Aa & \Bb
\end{pmatrix} T \vec{c} = 0 \big\} \ \ \text{is Lag.\ in\ } (\C^{2q},
\omega_E)\\
& \Longleftrightarrow \quad  \big\{ \vec{c} \in \C^{2 q} \, |\,
\begin{pmatrix} \Aa' & \Bb
\end{pmatrix} \vec{c} = 0 \big\} \ \ \text{is Lag.\ in\ } (\C^{2q},
\omega_E),
\end{align*}
which holds,  by Proposition \ref{p:relation}, if and only if
$\Aa' \Bb^*$ is self-adjoint.
\end{proof}

As seen in the formula \eqref{phiasymp1} of Proposition
\ref{p:maxdom}, the $\la_\ell = -\frac14$ eigenvalues of
$A_\Gamma$ and the $- \frac14 < \la_\ell < \frac34$ eigenvalues of
$A_\Gamma$ give rise to rather distinct components of
$\mathrm{dom}_{\max}(\Delta)$. For this reason, it is natural to
separate Lagrangian subspaces of $V$ into $\la_\ell = -\frac14$
components and $- \frac14 < \la_\ell < \frac34$ components.  With
this discussion in mind, we call a Lagrangian subspace $L \subset
V$ \emph{decomposable} if $L = L_0 \oplus L_1$ where $L_0$ is an
arbitrary Lagrangian subspace of $\bigoplus_{\la_\ell = -\frac14}
E_\ell^+ \oplus E_\ell^-$ and $L_1$ is an arbitrary Lagrangian
subspace of $\bigoplus_{-\frac14<\la_\ell<\frac34} E_\ell^+ \oplus
E_\ell^-$.

The characterization of all such $L_0,L_1$ follows from
Proposition \ref{p:relation}.

\begin{corollary} \label{cor-L01}
The components $L_0$ and $L_1$ of a decomposable Lagrangian $L =
L_0 \oplus L_1 \subset V$ can be characterized by matrices
$\begin{pmatrix} \Aa_0 & \Bb_0
\end{pmatrix}$ (with $\Aa_0$ and $\Bb_0$ $q_0 \times q_0$ matrices)
and $\begin{pmatrix} \Aa_1 & \Bb_1 \end{pmatrix}$ (with $\Aa_1$
and $\Bb_1$ $q_1 \times q_1$ matrices) via
\[
L_0 \cong \left\{ \vec{c} \in \C^{2 q_0} \, |\, \begin{pmatrix}
\Aa_0 & \Bb_0
\end{pmatrix} \vec{c} = 0 \right\} \subset \C^{2 q_0},
\]
and
\[
L_1 \cong \left\{ \vec{c} \in \C^{2 q_1} \, |\, \begin{pmatrix}
\Aa_1 & \Bb_1
\end{pmatrix} \vec{c} = 0 \right\} \subset \C^{2 q_1}
\]
where the matrix $\begin{pmatrix} \Aa_0 & \Bb_0 \end{pmatrix}$ has
rank $q_0$ and $\Aa_0 \Bb_0^* = \Bb_0 \Aa_0^*$, and
$\begin{pmatrix} \Aa_1 & \Bb_1 \end{pmatrix}$ has rank $q_1$ and
$\Aa_1 \Bb_1^* = \Bb_1 \Aa_1^*$.
\end{corollary}

In the Introduction we discussed split-type Lagrangians. Here, we
say that the Lagrangian $L_0 \subset \bigoplus_{\la_\ell =
-\frac14} E_\ell^+ \oplus E_\ell^-$ is of \emph{split-type} if it
can be written as $\bigoplus_{\la_\ell = - \frac14} L_{\ell}$ with
$L_\ell$ a Lagrangian subspace of $E_\ell^+ \oplus E_\ell^-$. In
the following proposition, we characterize all such Lagrangians
$L_\ell$.

\begin{proposition} \label{p:Ltheta}
$L \subset \C^2$ is Lagrangian if and only if $L = L_\theta$ for
some $\theta \in [0,\pi)$ where
\[
L_\theta = \{ (x,y) \in \C^2\ |\ \cos \theta\, x + \sin \theta\, y
= 0 \}.
\]
\end{proposition}
\begin{proof}
According to Proposition \ref{p:relation}, we have $L =  \left\{
(x,y) \in \C^{2} \, |\, a\, x + b\, y = 0 \right\} \subset
\C^{2}$, where $a,b \in \C$ with $a\, \overline{b} =
\overline{a}\, b$; that is, $a \, \overline{b} \in \R$. If $b =
0$, then $L_\ell \cong \left\{ (x,y) \in \C^{2} \, |\, 1 \cdot x +
0 \cdot y = 0 \right\}$, and we can take $\theta = 0$. If $b \ne
0$, then multiplying $a\, x + b\, y = 0$ by $\overline{b}$ we can
write
\[
L = \left\{ (x,y) \in \C^{2} \, |\, a\, \overline{b}\, x + |b|^2\,
y = 0 \right\} = \left\{ (x,y) \in \C^{2} \, |\, \alpha\, x +
\beta\, y = 0 \right\} ,
\]
where
\[
\alpha = \frac{a \, \overline{b}}{\sqrt{(a \, \overline{b})^2 +
(|b|^2)^2}} \ , \quad \beta = \frac{|b|^2}{\sqrt{(a \,
\overline{b})^2 + (|b|^2)^2}}  .
\]
Note that $\alpha,\beta \in \R$ and $\alpha^2 + \beta^2 = 1$. It
follows that $(\alpha,\beta) = (\cos \theta, \sin \theta)$ for
some $\theta \in \R$. If $\theta \in (\pi, 2 \pi]$, then we can
replace $(\alpha,\beta)$ by $(-\alpha,-\beta)$ to ensure that
$\theta \in [0,\pi)$. This completes our proof.
\end{proof}

\section{The model problems}
\label{s:modelproblems}
For the rest of this paper, unless stated otherwise, \emph{we fix
an arbitrary Lagrangian $L$ in $V$.} In this section we analyze
the eigenvalue equation for the model problem.

\subsection{The model operator}\label{subs:model}

In the last section we saw that only the  eigenvalues of
$A_\Gamma$ in the interval $[- \frac14, \frac34)$ are involved in
the various self-adjoint extensions of $\Delta$. For this reason,
in this section as a first step to prove our main results we shall
analyze the projection of
\[
\Delta|_{\mathcal{U}} = -\frac{d^2}{dr^2}+ \frac{1}{r^2} A_\Gamma
\]
onto the eigenspaces of $A_\Gamma$ with eigenvalues in
$[-\frac14,\frac34)$. Recall from Corollary \ref{cor-L}, with $q =
q_0 + q_1$, that the Lagrangian $L$ can be identified with the
null space of a $q \times 2q$ matrix $\begin{pmatrix} \Aa & \Bb
\end{pmatrix}$ of full rank with $\Aa$ and $\Bb$ $q \times q$
matrices such that $\Aa' \Bb^*$ is self-adjoint. Then writing
$A_\Gamma$ as a diagonal matrix with respect to its eigenfunctions
with eigenvalues in $[- \frac14 , \frac34)$, we shall consider the
operator
\[
\Ll := -\frac{d^2}{dr^2} + \frac{1}{r^2} A \qquad\text{over}\quad
[0,R],
\]
where $R > 0$ is arbitrary, but fixed, and $A$ is the $q \times q$
matrix
\[
A = \begin{pmatrix} - \frac14 \Id_{q_0}& 0 \\
0 & \begin{array}{ccccc}
\la_{q_0 + 1} & 0 & 0 & \cdots & 0\\
0 & \la_{q_0 + 2} & 0 & \cdots & 0\\
0 & 0 & \la_{q_0 + 3} & \cdots & 0 \\
0 & 0 & 0 & \ddots & 0\\
0 & 0 & 0 & \cdots & \la_{q_0 + q_1}
\end{array}
\end{pmatrix}.
\]
We put Dirichlet conditions at the right end $r = R$ of the
interval $[0,R]$. Then according to Proposition \ref{p:maxdom}, we
have

\begin{proposition} \label{p:maxdommodel}
$\phi \in \mathfrak{D}_\max$, the maximal domain of $\Ll$, if and
only if $\phi(R)=0$ and $\phi$ has the following form:
\begin{align} \label{phiasymp1model}
\phi =& \sum_{\ell = 1}^{q_0} \left\{ {c_\ell(\phi)} \, r^{
\frac12} e_\ell + {c_{q + \ell}(\phi)} \,
r^{ \frac12}\log r\, e_\ell \right\}\\
&\qquad + \sum_{\ell = 1}^{q_1} \left\{ c_{q_0 + \ell}(\phi) \,
r^{\nu_\ell + \frac12} e_{q_0 + \ell} + c_{q + q_0 + \ell}(\phi)
\, r^{-\nu_\ell + \frac12} e_{q_0 + \ell} \right\} +
\widetilde{\phi} ,\notag
\end{align}
where $\nu_j := \sqrt{\la_{q_0 + j} + \frac14}\ >\ 0$, $e_\ell$ is
the column vector with $1$ in the $\ell$-th slot and $0$'s
elsewhere, the $c_j(\phi)$'s are constants, and the
$\widetilde{\phi}$ is continuously differentiable on $[0,R]$ such
that $\widetilde{\phi}(r)=\mathcal{O}(r^{\frac32})$ and
$\widetilde{\phi}'(r) = \mathcal{O}(r^{\frac12})$ near $r=0$, and
$\Ll \widetilde{\phi} \in L^2([0,R],\C^q)$.
\end{proposition}

We dropped the factors $\frac{1}{\sqrt{2 \nu_\ell}}$, which appear
in the statement of Proposition \ref{p:maxdom}, from the terms in
\eqref{phiasymp1model} for $\phi$. Then as a consequence of
Theorem \ref{thm-sae}, we know that
\[
\Ll_L : \mathfrak{D}_L \to L^2([0,R], \C^{q})\ \, \text{is
self-adjoint,}\ \, \mathfrak{D}_L = \{\, \phi \in
\mathfrak{D}_\max\, |\, \vec{\phi} \in L \, \},
\]
where $\phi$ has the form in \eqref{phiasymp1model} with
$\vec{\phi} = (c_1(\phi),c_2(\phi),\ldots, c_{2q}(\phi))^t$. In
terms of the matrices $\Aa$ and $\Bb$, we can also write
\begin{equation} \label{DAB}
\mathfrak{D}_L = \{ \, \phi  \in \mathfrak{D}_\max\, |\,
\begin{pmatrix} \Aa & \Bb
\end{pmatrix} \vec{\phi} = 0 \ \} .
\end{equation}

\subsection{Eigenvalue equation}\label{subs:eigeneqn}

To analyze the $\zeta$-function of $\Ll_L$, we  derive an equation
for the eigenvalues of $\Ll_L$. For this, we first find solutions
to the equation
\[
(\Ll_L - \mu^2) \phi = 0 .
\]
As the reader can easily check, this is just a system of Bessel
equations as described in \cite[p.\ 362]{BAbM-StI92}, whose
solution (after judiciously choosing the constants for later
convenience) can be taken to be of the form
\begin{multline} \label{phi}
\phi = \sum_{\ell = 1}^{q_0} \Big\{ c_\ell(\phi) \, r^{\frac12} \,
J_{0}(\mu r) \, e_\ell + c_{q + \ell}(\phi) \, r^{\frac12} \,
\widetilde{J}_0(\mu r) \, e_\ell \Big\}\\+ \sum_{\ell = 1}^{q_1}
\Big\{ 2^{\nu_\ell} \Gamma(1 + \nu_\ell) c_{q_0 + \ell}(\phi) \,
\mu^{-\nu_\ell} \, r^{\frac12} \, J_{\nu_\ell}(\mu r) \, e_{q_0 + \ell} \\
+ 2^{-\nu_\ell} \Gamma(1 - \nu_\ell) c_{q + q_0 + \ell}(\phi) \,
\mu^{\nu_\ell} \, r^{\frac12} \,J_{-\nu_\ell}(\mu r)\, e_{q_0 +
\ell} \Big\} ,
\end{multline}
where $J_v(z)$ denotes the Bessel function of the first kind and
\begin{equation}\label{def-Jtilde}
\widetilde{J}_0(\mu r) := \frac{\pi}{2}  Y_0(\mu r) - (\log \mu -
\log 2 + \gamma)\, J_0(\mu r) ,
\end{equation}
with $Y_0(z)$ the Bessel function of the second kind. For
notational convenience let us introduce $J_{+0} (\mu R ) = J_0
(\mu R)$ and $J_{-0} (\mu R) = \widetilde{J}_0 (\mu R)$.

Define $q \times q$ matrices $J_+(\mu), J_-(\mu)$ by \tiny
\[
J_\pm(\mu) := \begin{pmatrix} J_{\pm 0}(\mu R) \Id_{q_0} & 0 & \cdots & 0 \\
0 & 2^{\pm\nu_1} \Gamma(1 \pm \nu_1) \, \mu^{\mp\nu_1}
J_{\pm\nu_1}(\mu R) & \cdots & 0 \\
0 &  0 & \cdots & 0\\
0 &  0 & \ddots & 0\\
0 &  0 & \cdots & 2^{\pm\nu_{q_1}} \Gamma(1 \pm \nu_{q_1}) \,
\mu^{\mp\nu_{q_1}} J_{\pm\nu_{q_1}}(\mu R) \end{pmatrix}.
\]
\normalsize
In the following proposition, we determine an
eigenvalue equation for the $\mu$'s.

\begin{proposition} \label{prop-Fmu}
$\mu^2$ is an eigenvalue of $\Ll_L$ if and only if
\[
F(\mu) := \det \begin{pmatrix} \Aa & \Bb \\
J_+(\mu) & J_-(\mu) \end{pmatrix} = 0.
\]
\end{proposition}
\begin{proof} Imposing the
Dirichlet condition at $r = R$ on $\phi$ of the form \eqref{phi},
we obtain
\[
c_\ell(\phi)\, J_0(\mu R) + c_{q + \ell}(\phi)\,
\widetilde{J}_0(\mu R) = 0 \ ,\quad \ell = 1,\ldots, q_0,
\]
and
\begin{multline*}
2^{\nu_\ell} \Gamma(1 + \nu_\ell) c_{q_0 + \ell}(\phi) \,
\mu^{-\nu_\ell} J_{\nu_\ell}(\mu R) \\ + 2^{-\nu_\ell} \Gamma(1 -
\nu_\ell) c_{q + q_0 + \ell}(\phi) \, \mu^{\nu_\ell}
J_{-\nu_\ell}(\mu R) = 0\ ,\quad \ell = 1,\ldots, q_1.
\end{multline*}
We can summarize these two equations as
\begin{equation}\label{def-eq1}
\begin{pmatrix} J_+(\mu) & J_-(\mu) \end{pmatrix} \vec{\phi} = 0,
\end{equation}
where $\vec{\phi} = (c_1(\phi), c_2(\phi), \ldots,
c_{2q}(\phi))^t$. Now recall that \cite[p.\ 360]{BAbM-StI92}
\begin{equation}\label{Jasymp0}
z^{-v} J_v(z) = \sum_{k = 0}^\infty \frac{(-1)^k z^{2k}}{2^{v + 2
k} k! \, \Gamma ( v + k + 1)}\sim \frac{1}{2^v \Gamma(1 + v)}
\hspace{-.2em} \left( 1 - \frac{z^2}{4 (1 + v)} + \cdots \right),
\end{equation}
and
\begin{equation}\label{Y0}
\frac{\pi}{2} Y_0(z) = \big( \log z - \log 2 + \gamma \big) J_0(z)
- \sum_{k = 1}^\infty \frac{H_k (- \frac14 z^2 )^k}{(k!)^2} ,
\end{equation}
where $H_k := 1 + \frac12 + \cdots + \frac{1}{k}$. Combining
\eqref{def-Jtilde}, \eqref{Jasymp0} with $v=0$, and \eqref{Y0}, we
get
\begin{equation} \label{J0tilde}
\widetilde{J}_0(\mu r) = (\log r)\, J_0(\mu r) -  \sum_{k =
1}^\infty \frac{H_k (- \frac14 (\mu r)^2 )^k}{(k!)^2} = \log r +
\mathcal{O}(r) .
\end{equation}
From \eqref{phi}, \eqref{Jasymp0} and \eqref{J0tilde}, it follows
that
\begin{multline*} \phi \, \sim\, \sum_{\ell = 1}^{q_0} \left\{ c_\ell(\phi) \,
r^{\frac12} e_\ell + c_{q + \ell}(\phi) \, r^{\frac12}\, \log r\,
e_\ell \right\}\\ + \sum_{\ell = 1}^{q_1} \left\{
c_{q_0+\ell}(\phi) \, r^{\nu_\ell + \frac12} e_{q_0+\ell} + c_{q
+q_0+ \ell}(\phi) \, r^{-\nu_\ell + \frac12} e_{q_0+\ell} \right\}
\quad \text{near}\ \ r=0 .
\end{multline*}
In particular, by \eqref{DAB}, $\phi$ in $\mathfrak{D}_L$
satisfies $\begin{pmatrix} \Aa & \Bb \end{pmatrix} \vec{\phi} =
0$, and therefore, in view of \eqref{def-eq1}, we conclude that
\[
\begin{pmatrix} \Aa & \Bb \\
J_+(\mu) & J_-(\mu) \end{pmatrix} \vec{\phi} = 0.
\]
For nontrivial $\vec{\phi}$, this equation can hold if and only if
the matrix in front of $\vec{\phi}$ is singular. This completes
our proof.
\end{proof}

\subsection{Asymptotics of $F(\mu)$}

In order to find relevant properties of the resolvent, the
heat-trace and the $\zeta$-function, we shall need the asymptotics
of $F(\mu)$ as $|\mu|\to\infty$.

\begin{proposition} \label{prop-asympF}
Let $\Upsilon \subset \C$ be a sector (closed angle) in the
right-half plane. Then as $|x| \to \infty$ with $x \in \Upsilon$,
we have
\begin{multline} \label{asymF}
F(i x) \sim (2 \pi R)^{-\frac q2} \prod_{j = 1}^{q_1} 2^{-\nu _{j}
} \Gamma (1 - \nu _{j} ) \, x^{|\nu| -\frac q2} \, e^{q x R}\,
(\widetilde{\gamma} - \log x)^{q_0} \times \\ p\Big( \big(
\widetilde{\gamma} - \log x \big)^{-1}, x^{-1} \Big) \, \Big( 1 +
\mathcal{O}(x^{-1}) \Big) ,
\end{multline}
where $\widetilde{\gamma} = \log 2 -\gamma$,
$p(x,y)$ is given in \eqref{polyp},
$\mathcal{O}(x^{-1})$ is a power series in $x^{-1}$, and
\begin{multline} \label{asymF2}
\frac{d}{d x} \log F(i x) \sim q R  + \frac{q_0 - j_0}{x (\log x - \widetilde{\gamma})}  \\
+\sum c_{\ell \xi}\, x^{-2 \xi - 1} \Big\{ \ell
(\widetilde{\gamma} - \log x)^{-\ell - 1} - 2 \xi
(\widetilde{\gamma} - \log x)^{- \ell} \Big\}  +
\mathcal{O}(x^{-1}) ,
\end{multline}
with the same meaning for $\mathcal{O}(x^{-1})$ and where the
$c_{\ell \xi}$'s are the constants in \eqref{clx}.
\end{proposition}
\begin{proof}
Using the identity $(i z)^{-v}\, J_v(i z) = z^{-v} I_v(z)$, where
$I_v(z)$ is the modified Bessel function of the first kind, we can
write $F(i x) = \det \begin{pmatrix} \Aa & \Bb\\
J_+(i x) & J_-(i x) \end{pmatrix}$ where (we use the notation
$I_{\pm 0} (xR) = J_{\pm 0} (ixR)$) \tiny
\[
J_\pm(i x) = \begin{pmatrix} I_{ \pm 0}(x R) \Id_{q_0} & 0 & \cdots & 0 \\
0 & 2^{\pm\nu_1} \Gamma(1 \pm \nu_1) \, x^{\mp\nu_1}
I_{\pm\nu_1}( x R) & \cdots & 0 \\
0 &  0 & \cdots & 0\\
0 &  0 & \ddots & 0\\
0 &  0 & \cdots & 2^{\pm\nu_{q_1}} \Gamma(1 \pm \nu_{q_1}) \,
x^{\mp\nu_{q_1}} I_{\pm\nu_{q_1}}(x R) \end{pmatrix}.
\]
\normalsize  Factoring out $2^{-\nu_j} \Gamma(1 - \nu_j) \,
x^{\nu_j}I_{-\nu_j}(x R)$ from
the $(q + q_0 + j)$-th row of the matrix $\begin{pmatrix} \Aa & \Bb \\
J_+(i x) & J_-(i x) \end{pmatrix}$ we obtain
\begin{multline} \label{Fixformula}
F(i x) \, = \, \rho \prod_{j=1}^{q_1} x^{\nu_j} I_{-\nu_j}(x R) \times \\
\det \begin{pmatrix} \Aa & \Bb \\
\begin{array}{cc} I_0(x R)\, \Id_{q_0} & 0 \\ 0 & A(x)
\end{array}
&
\begin{array}{cc}
\widetilde{J}_0(i x R)\, \Id_{q_0} & 0 \\
0 & \Id_{q_1}
\end{array}  \end{pmatrix},
\end{multline}
where $\rho = \prod_{j=1}^{q_1} 2^{-\nu_j} \Gamma(1 - \nu_j)$ and
\tiny
\[
A (x) = \begin{pmatrix} \tau_1\, x^{-2\nu_1} \frac{ I_{\nu_1}(x
R)}{I_{-\nu_1}(x R)}
& 0 & \cdots & 0 \\
0 & \tau_2 \,x^{-2\nu_2}
\frac{ I_{\nu_2}(x R)}{I_{-\nu_2}(x R)} & \cdots & 0\\
0 & 0 & \ddots & 0\\
0 & 0 & \cdots & \tau_{q_1} \, x^{-2\nu_{q_1}}\frac{
I_{\nu_{q_1}}(x R)}{I_{-\nu_{q_1}}(x R)}
\end{pmatrix},
\]
\normalsize with $\tau_j = 2^{2\nu_j}\frac{ \Gamma(1 + \nu_j)}{
\Gamma(1 - \nu_j)}$. In order to find the asymptotics of $F(ix)$
in \eqref{Fixformula}, we shall determine the asymptotics of
$A(x)$ and then of $\widetilde{J}_0(x)$. To determine the
asymptotics of $A(x)$, we recall (see \cite[p.\ 377]{BAbM-StI92})
that as $z\to\infty$ with $z \in \Upsilon$, we have
\begin{equation} \label{Ivasymp}
I_{v}(z) \sim \frac{e^z}{\sqrt{2\pi z}}\Big(1-\frac{4v^2-1}{8z} +
\mathcal{O}(z^{-2})\Big) \sim \frac{e^z}{\sqrt{2\pi z}} \Big(1 +
\mathcal{O}\Big(\frac{1}{z}\Big)\Big) ,
\end{equation}
where $\mathcal{O}\Big(\frac{1}{z}\Big)$ is a power series in
$\frac 1z$ and where only $v^2$'s occur in
$\mathcal{O}\Big(\frac{1}{z}\Big)$. In particular, as $x\to\infty$
with $x \in \Upsilon$, we have $\frac{I_{\nu_j}(x R)}{I_{-\nu_j}(x
R)} \sim 1$. Therefore,
\begin{equation} \label{Aasymp}
A (x) \sim \begin{pmatrix} \tau_1\, x^{-2\nu_1}
& 0 & \cdots & 0 \\
0 & \tau_2 \, x^{-2\nu_2} & \cdots & 0\\
0 & 0 & \ddots & 0\\
0 & 0 & \cdots & \tau_{q_1} \, x^{-2\nu_{q_1}}
\end{pmatrix}.
\end{equation}
To determine the asymptotics of $\widetilde{J}_0(i x)$, note that
$J_0(i z) = I_0(z)$ and
\begin{align*}
\frac{\pi}{2} Y_0(i z) & = \big( \log (i z) - \log 2 + \gamma
\big) J_0(i z) - \sum_{k = 1}^\infty \frac{H_k (- \frac14 (i z)^2
)^k}{(k!)^2}\\
& = \big( \log z + i \frac{\pi}{2} - \log 2 + \gamma \big) I_0(z)
- \sum_{k = 1}^\infty \frac{H_k (\frac14 z^2 )^k}{(k!)^2} = i
\frac{\pi}{2} I_0(z) - K_0(z)  ,
\end{align*}
where
\[
K_0(z) := - \big( \log z - \log 2 + \gamma \big) I_0(z) + \sum_{k
= 1}^\infty \frac{H_k (\frac14 z^2 )^k}{(k!)^2}
\]
is the modified Bessel function of the second kind. Thus, we can
write
\begin{align*}
\widetilde{J}_{0}(i x R) &= \frac{\pi}{2}  Y_0(i x R) - (\log (i
x) - \log 2 + \gamma)\, J_0(i x R) \\ &= i \frac{\pi}{2}  I_0(x R)
- K_0(x R)   - \Big( i \frac{\pi}{2} + \log x - \log 2 + \gamma
\Big)\, I_0(x R) \\ &= - (\log x - \widetilde{\gamma}) I_0(x R) -
K_0(x R).
\end{align*}
By \cite[p.\ 378]{BAbM-StI92}, $K_0(x)$ is exponentially decaying
as $x \to \infty$ in $\Upsilon$, so
\[
\widetilde{J}_{0}(i x R)  = - (\log x - \widetilde{\gamma}) I_0(x
R) - K_0(x R) \sim (\widetilde{\gamma} - \log x) I_0(x R).
\]
Summarizing our work so far, we see from \eqref{Fixformula} that
\begin{multline*}
F(i x) \, \sim \, \rho \, x^{|\nu|} \prod_{j=1}^{q_1} I_{-\nu_j}(x R) \times \\
\det \begin{pmatrix} \Aa & \Bb \\
\begin{array}{cc} I_0(x R)\, \Id_{q_0} & 0 \\ 0 & A(x)
\end{array}
&
\begin{array}{cc}
(\widetilde{\gamma} - \log x) I_0(x R) \, \Id_{q_0} & 0 \\
0 & \Id_{q_1}
\end{array}  \end{pmatrix},
\end{multline*}
where $\rho = \prod_{j=1}^{q_1} 2^{-\nu_j} \Gamma(1 - \nu_j)$,
$|\nu| = \nu_1 + \cdots + \nu_{q_1}$, and $A(x)$ satisfies
\eqref{Aasymp}. Now factoring out $(\widetilde{\gamma} - \log x)
I_0(x R)$ from $F(ix)$ and using the definition of $p(x,y)$ in
\eqref{polyp} (with ``$x$" replaced with $(\widetilde{\gamma} -
\log x)^{-1}$ and ``$y$" replaced with $x^{-1}$), we obtain
\[
F(i x) \, \sim \, \rho \, x^{|\nu|} \prod_{j=1}^{q_1} I_{-\nu_j}(x
R) \, I_0(x R)^{q_0} \, (\widetilde{\gamma} - \log x)^{q_0} \,
p\Big( \big( \widetilde{\gamma} - \log x \big)^{-1}, x^{-1} \Big).
\]
In view of the asymptotics \eqref{Ivasymp} for $I_v(z)$, we get
\begin{multline*}
F(i x) \, \sim \, \rho \, x^{|\nu|}  \bigg( (2 \pi)^{-\frac12}
(x R)^{-\frac12} e^{x R} \bigg) ^{q_0+q_1} \, (\widetilde{\gamma} - \log x)^{q_0} \times \\
p\Big( \big( \widetilde{\gamma} - \log x \big)^{-1}, x^{-1} \Big)
\Big( 1 + \mathcal{O}(x^{-1}) \Big) ,
\end{multline*}
which is equivalent to
\begin{multline*}
F(i x) \sim (2 \pi R)^{-\frac q2} \prod_{j = 1}^{q_1} 2^{-\nu _{j}
} \Gamma (1 - \nu _{j} ) \, x^{|\nu| -\frac q2} \, e^{q x R}\,
(\widetilde{\gamma} - \log x)^{q_0} \times \\ p\Big( \big(
\widetilde{\gamma} - \log x \big)^{-1}, x^{-1} \Big) \, \Big( 1 +
\mathcal{O}(x^{-1}) \Big) ,
\end{multline*}
and the proof of our first asymptotic formula is complete. To
prove our second formula, recall from \eqref{pxy} that
\[
p(x,y) = a_{j_0 \alpha_0}\, x^{j_0}\, y^{2 \alpha_0} \Big( 1 +
\sum  b_{k \beta} \, x^k \, y^{2 \beta} \Big) ,
\]
so that (with ``$x$" replaced with $(\widetilde{\gamma} - \log
x)^{-1}$ and ``$y$" replaced with $x^{-1}$)
\begin{multline*}
F(i x) \sim \text{const.} \times x^{|\nu| -\frac q2 - 2
\alpha_0}\, e^{q x R}\,  (\widetilde{\gamma} - \log x)^{q_0 - j_0}
\times
\\ \Big( 1 +
\sum  b_{k \beta} \, (\widetilde{\gamma} - \log x)^{-k} \, x^{- 2
\beta} \Big) \, \Big( 1 + \mathcal{O}(x^{-1}) \Big) .
\end{multline*}
As in  \eqref{clx}, $\log \Big( 1 + \sum  b_{k \beta} x^k y^{2
\beta} \Big) = \sum c_{\ell \xi} \, x^\ell \, y^{2 \xi}$ so taking
the logarithm of $F(ix)$ we see that
\begin{multline*}
\log F(i x) \sim \text{const.} + q x R + \Big( |\nu| -\frac q2 - 2
\alpha_0 \Big) \log x + \big( q_0 - j_0 \big) \log
(\widetilde{\gamma} - \log x) +
\\ \sum c_{\ell \xi} \, (\widetilde{\gamma} - \log x)^{-\ell} \, x^{- 2 \xi} + \mathcal{O}(x^{-1})  ,
\end{multline*}
and taking the derivative of both sides completes our proof.
\end{proof}

\subsection{The log terms only case} \label{ssec-logonly}

Suppose that $q_1 = 0$ so that the only eigenvalues of $A_\Gamma$
in the critical interval $[-\frac14,\frac34)$ are the $-\frac14$
eigenvalues. In this case, we shall denote $\Aa$ by $\Aa_0$ and
$\Bb$ by $\Bb_0$ so that
\[
F(\mu) = \det \begin{pmatrix} \Aa_0 & \Bb_0 \\
J_{0}(\mu R)\,  \Id_{q_0} & \widetilde{J}_{0}(\mu R)\,  \Id_{q_0}
\end{pmatrix} .
\]
Recall from Section \ref{ssec-special} (see \eqref{p0(z)}) the
polynomial
\[
p_0(z) := \det \begin{pmatrix} \Aa_0 & \Bb_0 \\
\Id_{q_0} & \big( \widetilde{\gamma} - z \big) \, \Id_{q_0}
\end{pmatrix},\quad \text{where}\ \ \widetilde{\gamma} = \log 2 -
\gamma,
\]
which is a polynomial in the complex variable $z \in \C$ of degree
at most $q_0$. Then we can write
\[
\frac{p_0'(z)}{p_0(z)} = \sum_{k = 1}^\infty \frac{\beta_k}{z^k} ,
\]
where the series is absolutely convergent for $|z|$ sufficiently
large and where $\beta_1 = \mathrm{deg}\, p_0$.  In the case that
$q_1 = 0$, Proposition \ref{prop-asympF} can be written as
follows.

\begin{proposition} \label{prop-asympF0}
Suppose that $q_1 = 0$ and let $\Upsilon \subset \C$ be a sector
(closed angle) in the right-half plane. Then as $|x| \to \infty$
with $x \in \Upsilon$, we have
\begin{equation} \label{asymFinfinity}
F(i x) \sim  (2 \pi x R)^{-\frac{q_0}{ 2}} \, e^{q_0 x R}\,
p_0(\log x) \,\Big( 1 + \mathcal{O}(x^{-1}) \Big) ,
\end{equation}
where $\mathcal{O}(x^{-1})$ is a power series in $x^{-1}$, and
\begin{align}\label{asymFlog}
\frac{d}{d x} \log F(i x)  & \sim  \sum_{k = 1}^\infty
\frac{\beta_k}{x(\log x)^k} + q_0 R - \frac{q_0}{2 x} +
\mathcal{O}(x^{-2}) ,
\end{align}
where $\mathcal{O}(x^{-2})$ is a power series in $x^{-1}$ starting
from $x^{-2}$.
\end{proposition}
\begin{proof}
A direct application of \eqref{asymF} in Proposition
\ref{prop-asympF} with $q_1 = 0$ gives
\[
F(i x) \sim (2 \pi x R)^{-\frac{q_0}{2}}\, e^{q_0 x R}\,
(\widetilde{\gamma} - \log x)^{q_0} p\Big( \big(
\widetilde{\gamma} - \log x \big)^{-1} \Big) \,  \Big( 1 +
\mathcal{O}(x^{-1}) \Big) ,
\]
where $\mathcal{O}(x^{-1})$ is a power series in $x^{-1}$ and
\[
p(x) := \det
\begin{pmatrix} \Aa_0 & \Bb_0 \\
x\, \Id_{q_0} & \Id_{q_0}
\end{pmatrix}.
\]
By definition of $p_0$, we have $(\widetilde{\gamma} - \log
x)^{q_0} p\Big( \big( \widetilde{\gamma} - \log x \big)^{-1} \Big)
= p_0( \log x)$. This proves \eqref{asymFinfinity} and then taking
the logarithmic derivative of \eqref{asymFinfinity} gives
\eqref{asymFlog}.
\end{proof}

\section{The zeta function for the model problems}
\label{s:modelzeta}

\begin{figure}
\centering \includegraphics{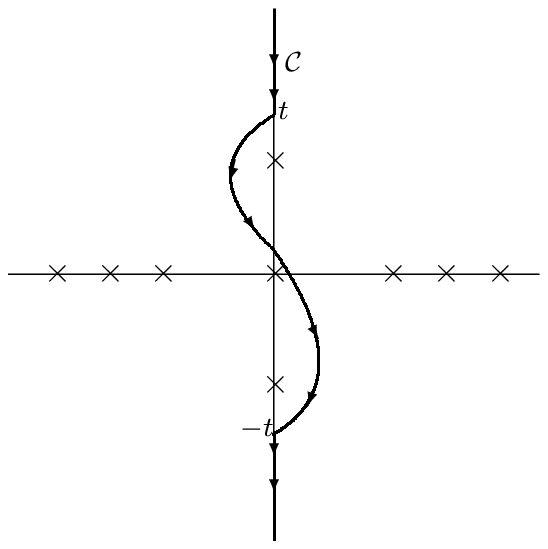} \caption{The
contour $\Cc$ for the zeta function. The $\times$'s represent the
zeros of $F(\mu)$. The squares of the $\times$'s on the imaginary
axis represent the negative eigenvalues of $\Ll_L$. Here, $t$ is
on the imaginary axis and $|t|^2$ is larger than the absolute
value of the negative eigenvalue of $\Ll_L$ (if one exists). The
contour $\Cc_t$ goes from $t$ to $-t$.} \label{fig-contstrange1}
\end{figure}

Working with the fixed Lagrangian $L$, we now analyze the zeta
function of $\Ll_L$, which by the Argument Principle \cite[p.\
123]{BCoJ78} is given by
\[
\zeta(s , \Ll_L ) = \frac{1}{2 \pi i} \int_\Cc \mu^{-2 s}
\frac{d}{d \mu} \log F(\mu) d \mu = \frac{1}{2 \pi i} \int_\Cc
\mu^{-2 s} \frac{F'(\mu)}{F(\mu)} d \mu,
\]
where $\Cc$ is a contour in the plane shown in Figure
\ref{fig-contstrange1}. Here we used that $\mu^2$ is an eigenvalue
of $\Ll_L$ if and only if $\mu$ is a zero of $F(\mu)$. By
Proposition \ref{prop-asympF}, the zeta function $\zeta(s,\Ll_L)$
is well-defined for $\Re s > \frac12$.

\subsection{A basic lemma}

In order to determine the exact structure of the analytic
continuation of $\zeta(s , \Ll_L )$, we need the following
fundamental result.

\begin{lemma} \label{lem-logint}
Let $c$ be a constant and let $|t|$ be sufficiently large so that
$\log x > c$ for $x \geq |t|$. Then for any $k \geq 0$ we can
write
\begin{equation} \label{intkpos}
\int_{|t|}^\infty x^{-2 s - 2 \xi - 1} \big(c - \log x)^k d x
\equiv \sum_{j = 0}^k \frac{\sigma_{k j}}{(s + \xi)^{j + 1}},
\end{equation}
modulo an entire function, where $\sigma_{k j} = (-1)^k
\binom{k}{j} \frac{j!}{2^{j+1}} (-c)^{k - j}$, and for any $k > 0$
we can write
\begin{equation} \label{intkneg}
\int_{|t|}^\infty x^{-2 s - 2 \xi - 1} \big(c - \log x)^{-k} d x
\equiv \frac{2^{k - 1}}{(k - 1)!}(s + \xi)^{k - 1} e^{-2 (s + \xi)
c} \log (s + \xi)
\end{equation}
modulo an entire function.
\end{lemma}
\begin{proof} Replacing $s$ by $s - \xi$, we can assume that $\xi = 0$
from the start. To analyze the first integral we first expand
$\big(c - \log x)^k$ using the binomial theorem:
\begin{align*}
\int_{|t|}^\infty x^{-2 s - 1} \big(c - \log x)^k d x  & = \sum_{j
= 0}^k (-1)^k \binom{k}{j} (-c)^{k - j} \int_{|t|}^\infty x^{-2 s
- 1} (\log x)^j d x .
\end{align*}
Thus, we are left to prove that
\begin{equation} \label{intlogxj}
\int_{|t|}^\infty x^{-2 s - 1} (\log x)^j d x \equiv
\frac{j!}{2^{j+1}} \cdot \frac{1}{s^{j+1}}
\end{equation}
modulo an entire function. However, since the integral
$\int_{1}^{|t|} x^{-2 s - 1} (\log x)^j d x$ is entire, we can
assume that the lower limit of the integral in \eqref{intlogxj} is
$1$. Now taking $j$ derivatives of both sides of the equality
$\int_{1}^\infty x^{-2 s - 1} d x = \frac{1}{2 \, s}$ with respect
to $s$, we obtain
\[
(-2)^j \int_{1}^\infty x^{-2 s - 1} (\log x)^j d x =
\frac{(-1)^{j} j!}{2 \, s^{j+1}},
\]
which proves \eqref{intlogxj}.

To prove the second claim in this proposition, we make the change
of variables $y = 2s (\log x - c)$ or $x = e^{c} e^{y/2s}$, and
obtain
\[
\int_{|t|}^\infty x^{-2 s - 1} \big(c - \log x \big)^{-k}\, d x =
(-1)^k \, e^{-2 s c} (2 s)^{k-1} \int_{2 s C}^\infty e^{- y}
\frac{d y}{y^k} .
\]
where $C := \log |t| - c$. Recall that the \emph{exponential
integral} is defined by (see \cite[p.\ 228]{BAbM-StI92} or
\cite[Sec.\ 8.2]{BGrI-RyI00})
\[
\Ei_k (z) := \int_{1}^\infty e^{-z u} \frac{d u}{u^k} = z^{k-1}
\int_{z}^\infty e^{-y} \frac{d y}{y^k}.
\]
Therefore,
\begin{align}
\int_{|t|}^\infty x^{-2 s - 1} \big(c - \log x \big)^{-k}\, d x =
\frac{(-1)^k \, e^{-2 s c}}{C^{k-1}} \, \Ei_k \big( 2 s C\big).
\label{intlogxk}
\end{align}
Also from \cite[p.\  229]{BAbM-StI92} or \cite[p.\
877]{BGrI-RyI00}, we have
\[
\Ei_k(z) = \frac{(-z)^{k-1}}{(k-1)!} \big\{ - \log z + \psi(k)
\big\}\ - \! \! \sum_{j = 0\, ,\, j \ne k - 1}^\infty
\frac{(-z)^{j-1}}{(j - k + 1) \, j!} ,
\]
where $\psi(1) := - \gamma$ and $\psi(k) : = - \gamma + \sum_{j =
1}^{k-1} \frac{1}{j}$ for $k > 1$. Hence,
\begin{align*} \label{Eikformula}
\Ei_k \big( 2 s C \big) = C^{k-1}\frac{(-2 s)^{k-1}}{(k-1)!}
\big\{ - \log (2 s C) + \psi(k) \big\} - \! \! \sum_{j = 0\, ,\, j
\ne k - 1}^\infty \frac{(-2 s C)^{j-1}}{(j - k + 1) \, j!}.
\end{align*}
Replacing this into \eqref{intlogxk} and simplifying, we obtain
\[
\int_{|t|}^\infty x^{-2 s - 1} \big(c - \log x \big)^{-k}\, d x
\equiv \frac{(2 s)^{k - 1}}{(k - 1)!} e^{-2 s c} \log s
\]
modulo an entire function. This completes our proof.
\end{proof}

\subsection{The $\z$-function}\label{sub-zetafnmodel}

We now prove the ``model problem version" of Theorem
\ref{thm-main} via the contour integration method
\cite{BKirK01,KM1,KM2}.

\begin{proposition} \label{prop-main}
Let $L \subset V$ be an \textbf{arbitrary} Lagrangian subspace of
$\C^{2q}$ and define $\mathscr{\pP}$ and $\mathscr{\lL}$ as in
\eqref{Spl} from the matrices $\Aa$ and $\Bb$ defining $L$. Then
the $\zeta$-function $\zeta(s,\Ll_L)$ extends from $\Re s
> \frac{1}{2}$ to a holomorphic function on $\C \setminus
(-\infty,0]$. Moreover, $\zeta(s,\Ll_L)$ can be written in the
form
\[
\zeta(s,\Ll_L) = \zeta_{\mathrm{reg}}(s,\Ll_L) +
\zeta_{\mathrm{sing}}(s,\Ll_L),
\]
where $\zeta_{\mathrm{reg}}(s,\Ll_L)$ has the ``regular" poles at
the ``usual" locations $s = \frac{1}{2}-k$ for $k \in \N_0$, and
where $\zeta_{\mathrm{sing}}(s,\Ll_L)$ has the following
expansion:
\begin{multline*}
\zeta_{\mathrm{sing}}(s,\Ll_L) = \frac{\sin (\pi s)}{\pi} \bigg\{
(j_0 - q_0) e^{-2 s (\log 2 - \gamma)} \log s + \sum_{\xi \in
\mathscr{\pP}} \frac{f_\xi(s)}{(s + \xi)^{|p_\xi| + 1}}\\
+ \sum_{\xi \in \mathscr{\lL}} g_\xi(s) \log (s + \xi) \bigg\},
\end{multline*}
where $j_0$ appears in \eqref{pxy} and $f_\xi(s)$ and $g_\xi(s)$
are entire functions of $s$ such that
\[
f_\xi(-\xi) = (-1)^{|p_\xi|+1} c_{p_\xi \xi} \, \xi\,
\frac{|p_\xi|!}{2^{|p_\xi|}}
\]
and
\[
g_\xi(s) = \begin{cases} c_{\ell_0, 0} \, \frac{2^{\ell_0}}{
(\ell_0 - 1)!} s^{\ell_0} + \mathcal{O}(s^{\ell_0 + 1}) & \text{if
\ \ $\xi = 0$,}\\ - c_{\ell_\xi \xi} \, \frac{\xi
2^{\ell_\xi}}{(\ell_\xi - 1)!} (s + \xi)^{\ell_\xi - 1} +
\mathcal{O}((s + \xi)^{\ell_\xi}) & \text{if \ \ $\xi > 0$,}
\end{cases}
\]
where the $c_{\ell \xi}$'s are the coefficients in \eqref{clx}.
\end{proposition}
\begin{proof}
With Figure \ref{fig-contstrange1} in mind, we write
\[
\int_\Cc = - \int_{t}^{0 + i \infty} + \int_{-t}^{0-i \infty} +
\int_{\Cc_t},
\]
where $\Cc_t$ is the curvy part of $\Cc$ from $t$ to $-t$, and
second, using that
\[
i^{-2s} = ( e^{i \pi/2} )^{-2 s} = e^{- i \pi s}\quad \text{and}
\quad (-i)^{-2s} = ( e^{-i \pi/2} )^{-2 s} = e^{i \pi s},
\]
we obtain the integral
\begin{align*}
\zeta(s,\Ll_L) & = \frac{1}{2 \pi i} \int_\Cc \mu^{-2 s}
\frac{d}{d\mu} \log F(\mu) \, d \mu\\
& = \frac{1}{2 \pi i} \bigg\{ - \int_{|t|}^\infty (i x)^{-2 s}
\frac{d}{d x} \log F(i x)\, d x + \int_{|t|}^\infty (-i x)^{-2 s}
\frac{d}{d x} \log F(-i x)\, d x \bigg\} \\
& \hspace{7.5cm} + \frac{1}{2 \pi i} \int_{\Cc_t} \mu^{-2 s}
\frac{F'(\mu)}{F(\mu)}\, d \mu\\ & = \frac{1}{2 \pi i} \Big( -
e^{- i \pi s} + e^{i \pi s}\Big) \int_{|t|}^\infty x^{-2 s}
\frac{d}{d x} \log F(i x)\, d x + \frac{1}{2 \pi i} \int_{\Cc_t}
\mu^{-2 s} \frac{F'(\mu)}{F(\mu)}\, d \mu ,
\end{align*}
or,
\begin{equation} \label{zetafn}
\zeta(s,\Ll_L) = \frac{\sin (\pi s)}{\pi} \int_{|t|}^\infty x^{-2
s} \frac{d}{d x} \log F(i x)\, d x + \frac{1}{2 \pi i}
\int_{\Cc_t} \mu^{-2 s} \frac{F'(\mu)}{F(\mu)}\, d \mu ,
\end{equation}
a formula that will be analyzed in a moment. The second integral
here is over a bounded contour so is an entire function of $s \in
\C$, so we are left to analyze the analytic properties of the
first integral in \eqref{zetafn}. To do so, recall the asymptotics
\eqref{asymF2} in Proposition \ref{prop-asympF}, which states that
for $x \to \infty$ we have
\begin{align}\label{e:G}
\frac{d}{d x} \log F(i x) &
\sim \frac{q_0 - j_0}{x (\log x - \widetilde{\gamma})} + G_1(x) +
G_2(x) + G_3(x) ,
\end{align}
where $\widetilde{\gamma} = \log 2 - \gamma$, $G_3(x)$ is a power
series in $x^{-1}$ starting with the constant term $q R$,
\[
G_1(x) := \sum_\xi \sum_{\ell \leq 0} c_{\ell \xi}\, x^{-2 \xi -
1} \Big\{ \ell (\widetilde{\gamma} - \log x)^{-\ell - 1} - 2 \xi
(\widetilde{\gamma} - \log x)^{- \ell} \Big\},
\]
and
\[
G_2(x) := \sum_\xi \sum_{\ell > 0} c_{\ell \xi}\, x^{-2 \xi - 1}
\Big\{ \ell (\widetilde{\gamma} - \log x)^{-\ell - 1} - 2 \xi
(\widetilde{\gamma} - \log x)^{- \ell} \Big\} .
\]
 Since
\[
\frac{\sin (\pi s)}{\pi} \int_{|t|}^\infty x^{-2 s - k}\, d x =
\frac{\sin (\pi s)}{\pi} \frac{x^{-2 s - k + 1}}{-2 s - k + 1}
\bigg|_{x = |t|}^\infty = \frac{\sin (\pi s)}{\pi} \frac{|t|^{-2 s
- k + 1}}{2 s + k - 1}
\]
which has poles at $s = \frac{1 - k}2$ for $s \notin \Z$, it
follows that
\begin{equation} \label{F3x}
\frac{\sin (\pi s)}{\pi} \int_{|t|}^\infty x^{-2 s} \frac{d}{d x}
\log G_3(i x)\, d x
\end{equation}
will contribute to the function $\zeta_{\mathrm{reg}}(s,\Ll_L)$ in
the statement of this proposition.
Setting $\xi = 0$ and $k = 1$ in equation \eqref{intkneg} in Lemma
\ref{lem-logint} we see that
\[
\int_{|t|}^\infty x^{-2 s} \frac{q_0 - j_0}{x (\log x -
\widetilde{\gamma})} d x \equiv - (q_0 - j_0) e^{-2 s
\widetilde{\gamma}} \log s
\]
modulo an entire function, which gives us the first term in
$\zeta_{\mathrm{sing}}(s,\Ll_L)$.

We now analyze $\int_{|t|}^\infty x^{-2 s} G_2(x)\, dx$. To do so,
we apply equation \eqref{intkneg} term-by-term to
\begin{multline*}
\int_{|t|}^\infty x^{-2 s} G_2(x)\, dx =\\ \int_{|t|}^\infty x^{-2
s} \sum_\xi \sum_{\ell > 0} c_{\ell \xi}\, x^{-2 \xi - 1} \Big\{
\ell (\widetilde{\gamma} - \log x)^{-\ell - 1} - 2 \xi
(\widetilde{\gamma} - \log x)^{- \ell} \Big\} dx,
\end{multline*}
and we see that, modulo an entire function,
\begin{align*}
\int_{|t|}^\infty & x^{-2 s} G_2(x)\, dx \equiv \sum_{\xi}
\sum_{\ell
> 0} c_{\ell \xi}\, \Big\{ \frac{\ell 2^\ell}{\ell!} \, e^{-2 (s +
\xi) \widetilde{\gamma}} \, (s + \xi)^\ell \log(s + \xi)\\ &
\qquad \qquad \qquad \qquad \qquad \qquad - \frac{\xi
2^{\ell}}{(\ell - 1)!} \, e^{-2 (s + \xi) \widetilde{\gamma}} \,
(s +
\xi)^{\ell - 1} \log(s + \xi) \Big\}\\
& = \sum_{\xi} \bigg( e^{-2 (s + \xi) \widetilde{\gamma}}
\sum_{\ell > 0} c_{\ell \xi} \Big\{\frac{\ell 2^\ell}{\ell!} \, (s
+ \xi)^\ell - \frac{\xi 2^{\ell}}{(\ell - 1)!} \, (s + \xi)^{\ell
- 1}
 \Big\}  \bigg) \log(s + \xi)  ,
\end{align*}
which can be written in the form $\sum_{\xi} g_\xi(s) \, \log (s +
\xi)$ where
\[
g_\xi(s) = e^{-2 (s + \xi) \widetilde{\gamma}} \sum_{\ell > 0}
c_{\ell \xi} \Big\{ \frac{2^\ell}{(\ell-1)!} \, (s + \xi)^\ell -
\frac{\xi 2^{\ell}}{(\ell - 1)!} \, (s + \xi)^{\ell - 1}
 \Big\} .
\]
From this explicit formula for $g_\xi(s)$, we see that
\[
g_\xi(s) = \begin{cases} c_{\ell_0, 0} \, \frac{2^{\ell_0}}{
(\ell_0 - 1)!} s^{\ell_0} + \mathcal{O}(s^{\ell_0 + 1}) & \text{if
$\xi = 0$,}\\ - c_{\ell_\xi \xi} \, \frac{\xi
2^{\ell_\xi}}{(\ell_\xi - 1)!} (s + \xi)^{\ell_\xi - 1} +
\mathcal{O}((s + \xi)^{\ell_\xi}) & \text{if $\xi > 0$,}
\end{cases}
\]
where we recall that $\ell_\xi := \min \{\ell > 0 \, |\, c_{\ell
\xi} \ne 0\}$.

We now analyze $\int_{|t|}^\infty x^{-2 s} G_1(x)\, dx$.
With $\sigma_{k j} = (-1)^k \binom{k}{j} \frac{j!}{2^{j+1}}
(-\widetilde{\gamma})^{k - j}$, from equation \eqref{intkpos} in
Lemma \ref{lem-logint} we can write, modulo an entire function,
\begin{multline}\label{e:G1}
\int_{|t|}^\infty x^{-2 s} \bigg( \sum_{\ell \leq 0} c_{\ell
\xi}\, x^{-2 \xi - 1} \Big\{ \ell (\widetilde{\gamma} - \log
x)^{-\ell - 1} - 2 \xi
(\widetilde{\gamma} - \log x)^{- \ell} \Big\} \bigg) dx\\
= \sum_{\ell \leq 0} c_{\ell \xi}\, \bigg\{ \sum_{j=0}^{|\ell| -1}
\frac{\ell \sigma_{|\ell| - 1, j}}{(s + \xi)^{j + 1}} - \sum_{j =
0}^{|\ell|} \frac{2\xi \sigma_{|\ell| j}}{(s + \xi)^{j + 1}}
\bigg\} = \frac{f_\xi(s)}{(s + \xi)^{|p_\xi| + 1}},
\end{multline}
where
\[
p_\xi := \min \{\ell \leq 0 \, |\, c_{\ell \xi} \ne 0\}\quad
\Longrightarrow \quad |p_\xi| = \max \{\, |\ell|\ |\ \ell \leq 0 \
\text{and}\ c_{\ell \xi} \ne 0\}
\]
and
\[
f_\xi(s) := \sum_{\ell \leq 0} c_{\ell \xi}\, \bigg\{
\sum_{j=0}^{|\ell| - 1} \ell \sigma_{|\ell| - 1, j} \, (s +
\xi)^{|p_\xi| - j} - \sum_{j = 0}^{|\ell|} 2 \xi \sigma_{|\ell| j}
\, (s + \xi)^{|p_\xi| - j} \bigg\}
\]
is entire. It follows that
\begin{multline*}
\int_{|t|}^\infty x^{-2 s} G_1(x)\, dx =\\ \int_{|t|}^\infty x^{-2
s} \sum_\xi \sum_{\ell \leq 0} c_{\ell \xi}\, x^{-2 \xi - 1}
\Big\{ \ell (\widetilde{\gamma} - \log x)^{-\ell - 1} - 2 \xi
(\widetilde{\gamma} - \log x)^{- \ell}
\Big\} dx \\
= \sum_\xi \frac{f_\xi(s)}{(s + \xi)^{|p_\xi| + 1}}.
\end{multline*}
Moreover, from the above explicit formula for $f_\xi(s)$
we see that
\[
f_\xi(-\xi) = -  2 c_{p_\xi \xi} \xi \sigma_{|p_\xi| , |p_\xi|} =
- 2 c_{p_\xi \xi} \xi (-1)^{|p_\xi|}
\frac{|p_\xi|!}{2^{|p_\xi|+1}} = (-1)^{|p_\xi|+1} c_{p_\xi \xi} \,
\xi\,
 \frac{|p_\xi|!}{2^{|p_\xi|}} .
\]
This completes our proof.
\end{proof}

\subsection{The decomposable case}

Suppose now that $L = L_0 \oplus L_1$ is decomposable where as in
Corollary \ref{cor-L01} $L_0$ is given by $q_0 \times q_0$
matrices $\Aa_0$ and $\Bb_0$ and $L_1$ is given by $q_1 \times
q_1$ matrices $\Aa_1$ and $\Bb_1$. Let us recall the polynomial
$p_0(z)$ introduced in \eqref{p0(z)} in Section \ref{ssec-special}
and consider the following result.

\begin{lemma} \label{lem-logint2} For $|t|$ sufficiently large so
that $p_0(\log x)$ has no zeros for $x \geq |t|$, we can write
\[
\int_{|t|}^\infty x^{-2 s} \frac{p_0'( \log x)}{x p_0(\log x)}\, d
x \equiv - f(s) \log s
\]
modulo an entire function, where $f(s)$ is the entire function
given explicitly by
\[
f(s) = \sum_{ k = 1 }^\infty \beta_k \frac{(- 2 s)^{k-1}}{(k-1)!}
\]
with the $\beta_k$'s the coefficients of the expansion of
$\frac{p_0'(z)}{p_0(z)} = \sum_{k = 1}^\infty \frac{\beta_k}{z^k}$
in \eqref{p0'p0}.
\end{lemma}
\begin{proof} Using the expansion $\frac{p_0'(z)}{p_0(z)} =
\sum_{k = 1}^\infty \frac{\beta_k}{z^k}$, we can write
\begin{equation} \label{intformula}
\int_{|t|}^\infty x^{-2 s} \frac{p_0'( \log x)}{x p_0(\log x)}\, d
x = \sum_{k = 1}^\infty \beta_k \int_{|t|}^\infty x^{-2 s}
\frac{1}{x \, (\log x)^k}\, d x .
\end{equation}
To analyze this integral we put $\xi = c = 0$ in formula
\eqref{intkneg} from Lemma \ref{lem-logint} to see
\[
\int_{|t|}^\infty x^{-2 s} \frac{1}{x(\log x)^{k}} d x \equiv -
\frac{(-2)^{k - 1}}{(k - 1)!} \, s^{k - 1} \log s = - \frac{(-2
s)^{k - 1}}{(k - 1)!} \, \log s,
\]
modulo an entire function. Replacing this formula into
\eqref{intformula} and simplifying, we obtain our result.
\end{proof}

Let us apply this theorem to the case when $q_0 = 1$. In this
case, by Proposition \ref{p:Ltheta} we have $\Aa_0 = \cos \theta$
and $\Bb_0 = \sin \theta$ for an angle $\theta \in [0,\pi)$,
therefore
\[
p_0(z) := \det \begin{pmatrix} \cos \theta & \sin \theta \\
1 & \big(\widetilde{\gamma} - z \big) \end{pmatrix} = - \cos
\theta \cdot z + \widetilde{\gamma} \cdot \cos \theta - \sin
\theta .
\]
Hence, with $\kappa:= \widetilde{\gamma} - \tan \theta  = \log 2 -
\gamma - \tan \theta$, we have
\begin{align}
\frac{p_0'(z)}{p_0(z)} = \frac{- \cos \theta}{- \cos \theta \cdot
z + \widetilde{\gamma} \cdot \cos \theta - \sin \theta} & =
\frac{1}{z - \widetilde{\gamma} + \tan \theta}   = \sum_{k =
1}^\infty \frac{\kappa^{k-1}}{z^k}. \label{q0=1case}
\end{align}
Thus, $\beta_k = \kappa^{k-1}$, so
\[
f(s) = \sum_{ k = 1 }^\infty \beta_k \frac{(- 2 s)^{k-1}}{(k-1)!}
= \sum_{ k = 1 }^\infty \frac{(- 2 s \kappa)^{k-1}}{(k-1)!} = e^{-
2 s \kappa} .
\]
Therefore, Lemma \ref{lem-logint2} reduces to

\begin{corollary} \label{cor-m=1caseint}
Suppose that $q_0 = 1$, $q_1 = 0$ and  $\theta \ne \frac{\pi}{2}$.
Then the coefficients $\beta_k$ in the expansion
$\frac{p_0'(z)}{p_0(z)}$ are given by $\beta_k = \kappa^{k-1}$. In
particular, for $|t|$ sufficiently large so that $p_0(\log x)$ has
no zeros for $x \geq |t|$, we can write
\[
\int_{|t|}^\infty x^{-2 s} \frac{p_0'( \log x)}{x p_0(\log x)}\, d
x \equiv - e^{-2 s \kappa} \log s
\]
modulo an entire function, where $\kappa = \log 2 - \gamma - \tan
\theta$.
\end{corollary}

From \eqref{p1} and \eqref{cxi} in Section \ref{ssec-special}, let
us recall that the polynomial $p_1(y)$ has the expression
\[
p_1(y) =  a_{\alpha_0}\, y^{2 \alpha_0} \Big( 1 + \sum b_\beta
y^{2 \beta} \Big),
\]
where the $\beta$'s are positive and
\begin{equation} \label{coeffcxi}
\log \Big( 1 + \sum  b_{\beta} y^{2 \beta} \Big) = \sum c_{\xi} \,
y^{2 \xi} ,
\end{equation}
and let $\mathscr{\pP} := \{ \xi \, |\, c_\xi \ne 0\}$. We now
prove the model problem versions of Theorems \ref{main thm} and
\ref{main thm2}.

\begin{proposition} \label{prop-zeta}
For an arbitrary decomposable Lagrangian $L \subset V$, the
$\zeta$-function $\zeta(s,\Ll_L)$ has the following form:
\[
\z(s,\Ll_L) = \z_{\mathrm{reg}}(s,\Ll_L) +
\zeta_{\mathrm{sing}}(s,\Ll_L)
\]
where $\zeta_{\mathrm{reg}}(s,\Ll)$ has the ``regular" poles at
the ``usual" locations $s = \frac{1}{2} -k$ for $k \in \N_0$, and
where $\zeta_{\mathrm{sing}}(s,\Ll_L)$ has the following
expansion:
\[
\zeta_{\mathrm{sing}}(s,\Ll_L) =  - \frac{ \sin (\pi s)}{\pi}f(s)
\log s \ +\ \frac{\sin (\pi s)}{\pi} \! \! \sum_{\xi \in
\mathscr{\pP}} \frac{f_\xi(s)}{s + \xi},
\]
where $f(s)$ is the entire function defined explicitly by $f(s) =
\sum_{ k = 1 }^\infty \beta_k \frac{(- 2 s)^{k-1}}{(k-1)!}$, and
the $f_\xi(s)$'s are entire functions such that $f_\xi(-\xi) = -
{c_{\xi} \,  \xi}$ with the $c_\xi$'s the coefficients in
\eqref{coeffcxi}.
\end{proposition}
\begin{proof}
Since $L = L_0 \oplus L_1$ is decomposable, it follows that
\[
\zeta(s, \Ll_L) = \zeta(s , \Ll_{L_0}) + \zeta(s, \Ll_{L_1}),
\]
where $\Ll_{L_0}$ is the operator $\Ll_L$ restricted to the
$-\frac14$ eigenspaces of $A$ and $\Ll_{L_1}$ is the operator
$\Ll_L$ restricted to the eigenspaces of $A$ in $(-\frac14,
\frac34)$. From \eqref{coeffcxi}, we can observe that $p_\xi=0$
for any $\xi$ and $\mathscr{\lL}=\emptyset$ for the operator
$\Ll_{L_1}$.  Hence, there are only terms $G_1(x)$ with $\ell=0$
and $G_3(x)$ in \eqref{e:G}, that is,
\begin{equation}\label{e:decom}
\frac{d}{d x} \log F(i x) \, \sim \, \sum_{k=0}^\infty b_k x^{-k}
+ \sum - 2 \xi c_{ \xi}\, x^{-2 \xi - 1}
\end{equation}
where $c_{\xi}$'s are the coefficients in \eqref{coeffcxi}. It
follows from the proof of Proposition \ref{prop-main}, in
particular, \eqref{e:G1} with $\ell=0$ that
\[
\z(s,\Ll_{L_1}) = \z_{\mathrm{reg}}(s,\Ll_{L_1}) +
\zeta_{\mathrm{sing}}(s,\Ll_{L_1})
\]
where $\zeta_{\mathrm{reg}}(s,\Ll_{L_1})$ has poles at $s =
\frac{1}{2}-k$ for $k \in \N_0$ and
\[
\zeta_{\mathrm{sing}}(s,\Ll_{L_1}) = \frac{\sin (\pi s)}{\pi} \!
\! \sum_{\xi \in \mathscr{\pP}} \frac{f_\xi(s)}{s + \xi},
\]
where the $f_\xi(s)$'s are entire functions of $s$ such that
$f_\xi(-\xi) = - {c_{\xi} \,  \xi}$.

Thus, it remains to analyze $\zeta(s,\Ll_{L_0})$. To do so, we
follow the proof of Proposition \ref{prop-main} up to equation
\eqref{zetafn}, for $|t| \gg 0$ we can write
\[
\zeta(s,\Ll_{L_0}) \equiv \frac{\sin (\pi s)}{\pi}
\int_{|t|}^\infty x^{-2 s} \frac{d}{d x} \log F_0(i x)\, d x
\]
modulo an entire function, where $F_0(\mu)$ is the function
$F(\mu)$ defined in Proposition \ref{prop-Fmu} in the case that
$q_1 = 0$, $\Aa = \Aa_0$, and $\Bb = \Bb_0$. By \eqref{asymFlog}
in Proposition \ref{prop-asympF0}, we have
\[
\frac{d}{d x} \log F_0(i x) \sim \frac{p_0'(\log x)}{x\, p_0(\log
x)} + G(x),
\]
where $G(x)$ is a power series in $x^{-1}$ starting with a
constant term. Just as we noticed in \eqref{F3x} for $G_3(x)$
(which has the same asymptotics as $G(x)$) in the proof of
Proposition \ref{prop-asympF0}, the integral $\frac{\sin (\pi
s)}{\pi} \int_{|t|}^\infty x^{-2 s} G(i x)\, d x$ will contribute
to the function $\zeta_{\mathrm{reg}}(s,\Ll_{L_0})$ in the
statement of this proposition. Finally, invoking Lemma
\ref{lem-logint2}:
\[
\int_{|t|}^\infty x^{-2 s} \frac{p_0'( \log x)}{x p_0(\log x)}\, d
x \equiv - f(s) \log s
\]
modulo an entire function, where $f(s) = \sum_{ k = 1 }^\infty
\beta_k \frac{(- 2 s)^{k-1}}{(k-1)!}$, completes the proof.
\end{proof}

Corollary \ref{cor-m=1caseint} implies

\begin{corollary}
\label{cor-zeta} Suppose that $q_0 = 1$, $q_1 = 0$ and $\theta \ne
\frac{\pi}{2}$. Then the zeta function $\zeta(s,\Delta_\theta)$
can be written in the form
\[
\zeta(s,\Delta_\theta) = - \frac{ \sin (\pi s)}{\pi}e^{- 2 s
\kappa} \log s + \zeta_\theta(s),
\]
where $\kappa = \log 2 - \gamma - \tan \theta$ and
$\zeta_\theta(s)$ extends from $\Re s > \frac12$ to a holomorphic
function on $\C$ with poles at $s = \frac{1}{2}-k$ for $k \in
\N_0$. In particular, $\zeta(s , \Delta_\theta)$ has $s = 0$ as a
logarithmic branch point. In the case that $\theta =\frac{ \pi}
2$, the $\zeta$-function $\zeta(s,\Delta_\theta)$ has the
properties of $\zeta_\theta(s)$.
\end{corollary}

The last statement for $\theta=\frac{\pi}{2}$ in Corollary
\ref{cor-zeta} follows from the results in \cite{FMP}.

\section{The resolvent and heat kernel for the model problems}
\label{s:model}

In this section we analyze the resolvent and heat kernel
expansions for the model problems, which will be of great use for
the general case.

\subsection{The resolvent}

\begin{figure}
\centering \includegraphics{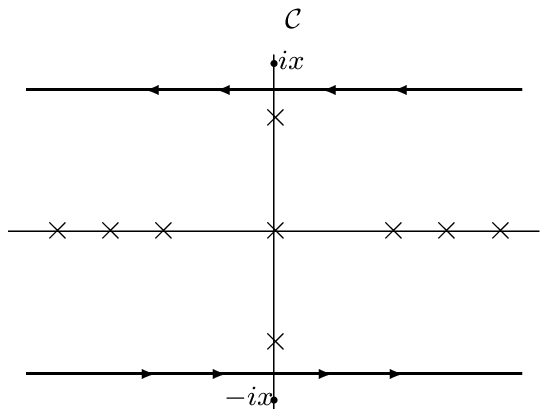} \caption{The new
contour $\Cc$.} \label{fig-contres}
\end{figure}

Using the new contour $\Cc$ shown in Figure \ref{fig-contres}, we
see that if $\{\mu_j^2\}$ denote the eigenvalues of $\Ll_L$, then
by an application of the Argument Principle, we have
\begin{align*}
2 \Tr( \Ll_L + x^2)^{-1} = 2 \sum_{j = 1}^\infty \frac{1}{\mu_j^2
+ x^2} &  = \frac{1}{2 \pi i} \int_\gamma (\mu^2 + x^2)^{-1}
\frac{d}{d \mu} \log F(\mu) \, d \mu ,
\end{align*}
where $|x|^2$ is larger than the absolute value of the negative
eigenvalues of $\Ll_L$ (if one exists). The factor of $2$ on the
left hand side is a result of all eigenvalues being enclosed
twice. Using this formula we can express the trace of the
resolvent in terms of $F(i x)$ in the following important and
remarkable theorem.

\begin{theorem} \label{thm-res}
We have
\[
2 x \Tr( \Ll_L + x^2)^{-1} = \frac{d}{d x} \log F(i x) .
\]
for all complex $x \in \C$ for which either (and hence both) sides
make sense.
\end{theorem}
\begin{proof} Deforming the contour as in Figure
\ref{fig-contres2}
\begin{figure}[b]
\centering \includegraphics{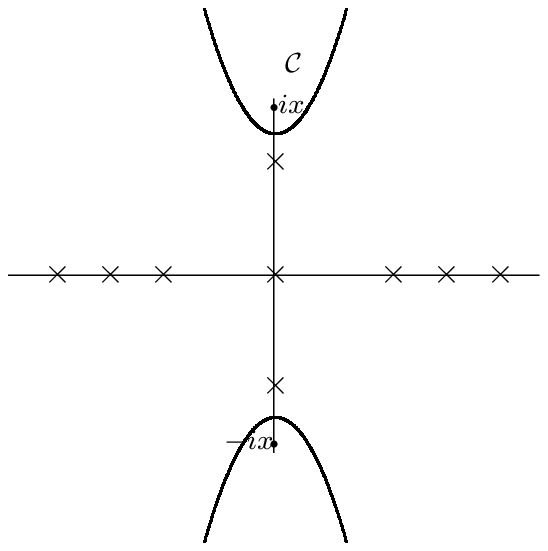} \caption{Deforming
the contour $\Cc$.} \label{fig-contres2}
\end{figure}
and using Cauchy's formula, we obtain
\begin{align*}
2 \Tr( \Ll_L + x^2)^{-1} & = \frac{1}{2 \pi i} \int_\gamma
\frac{1}{(\mu - i x)(\mu + i x)} \frac{F'(\mu)}{F(\mu)} d \mu\\ &
= - \frac{1}{2 i x} \frac{F'(ix)}{F(ix)} - \frac{1}{-2 i x}
\frac{F'(-ix)}{F(-ix)}  = \frac{i}{ x} \frac{F'(ix)}{F(ix)} =
\frac{1}{x} \frac{d}{d x} \log F(ix),
\end{align*}
where we used the fact that $F(\mu)$ is an even function of $\mu$.
Indeed, to see this observe that, by definition, $F(\mu)$ is
expressed in terms of $\mu^v J_{-v}(\mu R)$ with appropriate $v$'s
and the function $\widetilde{J}_0(\mu R)$, which are even
functions by \eqref{def-Jtilde}, \eqref{Jasymp0} and \eqref{Y0}.
This proves that $2 x \Tr( \Ll_L + x^2)^{-1} = \frac{d}{d x} \log
F(i x)$ at least when $x$ is real and $x \gg 0$. However, by
analytic continuation, both sides must still be equal for all
complex $x$ for which both sides are defined.
\end{proof}
\begin{remark} \em
This theorem is really quite remarkable because it tells us how to
\emph{immediately} evaluate traces of resolvents from simply
knowing an implicit eigenvalue equation! There are many
applications of this theorem that will appear elsewhere.
\end{remark}

Using this theorem, we can now prove

\begin{proposition}\label{prop-asympresol}
Let $L \subset V$ be an \textbf{arbitrary} Lagrangian subspace of
$\C^{2q}$ and let $\Lambda \subset \C$ be any sector (solid angle)
not intersecting the positive real axis. Then as $|\la| \to
\infty$ with $\la \in \Lambda$, we have
\begin{multline*}
\Tr( \Ll_L - \la)^{-1}\,  \sim \, \sum_{k = 1}^\infty a_k
(-\la)^{-\frac k2} + \frac{q_0 - j_0}{(-\la) (\log (-\la) - 2
\widetilde{\gamma})} \\ - \frac{d}{d \la} \left\{ \sum 2^\ell
c_{\ell \xi}\, (-\la)^{- \xi} \Big(2 \widetilde{\gamma} - \log
(-\la) \Big)^{-\ell} \right\},
\end{multline*}
where the $a_k$ coefficients are independent of $L$ and the
$c_{\ell \xi}$'s are given in \eqref{clx}.
\end{proposition}
\begin{proof}
By Proposition \ref{prop-asympF} (see equation \eqref{asymF2}), we
have
\begin{multline*}
\frac{d}{d x} \log F(i x) \,  \sim \, \sum_{k = 0}^\infty b_k
x^{-k}
+ \frac{q_0 - j_0}{x (\log x - \widetilde{\gamma})} \  \\
+\sum c_{\ell \xi}\, x^{-2 \xi - 1} \Big\{ \ell
(\widetilde{\gamma} - \log x)^{-\ell - 1} - 2 \xi
(\widetilde{\gamma} - \log x)^{- \ell} \Big\}
\end{multline*}
for some coefficients $b_k$ that, by the proof of Proposition
\ref{prop-asympF}, are independent of $L$. Therefore, by Theorem
\ref{thm-res}, we have
\begin{multline*}
2 x \Tr( \Ll_L + x^2)^{-1}  \,  \sim \, \sum_{k = 0}^\infty
b_k x^{-k} + \frac{q_0 - j_0}{x (\log x - \widetilde{\gamma})}  \\
+\sum c_{\ell \xi}\, x^{-2 \xi - 1} \Big\{ \ell
(\widetilde{\gamma} - \log x)^{-\ell - 1} - 2 \xi
(\widetilde{\gamma} - \log x)^{- \ell} \Big\}  .
\end{multline*}
Dividing by $2x$ and then setting $x = (-\la)^{\frac12}$, with
$a_k = (1/2) b_{k-1}$ we obtain
\begin{align*}
\Tr( \Ll_L - \la)^{-1}  \, & \sim \,
\sum_{k = 1}^\infty a_k (-\la)^{-\frac k2} + \frac{q_0 -
j_0}{(-\la) (\log (-\la) - 2\widetilde{\gamma})} \\ &
\hspace{3.5cm} -\frac{d}{d \la} \left\{ \sum c_{\ell \xi}\,
(-\la)^{- \xi} \Big(\widetilde{\gamma} - \frac12 \log (-\la)
\Big)^{-\ell} \right\},
\end{align*}
which is equivalent to our desired result.
\end{proof}

For decomposable $L$, we have

\begin{proposition} \label{prop-restracede}
Let $\Lambda \subset \C$ be any sector (solid angle) not
intersecting the positive real axis. Then for an arbitrary
decomposable Lagrangian $L$, as $|\la| \to \infty$ with $\la \in
\Lambda$ we have
\begin{align*}
\Tr (\Ll_L - \lambda)^{-1} \sim & \sum_{k=1}^\infty a_k \, (-
\lambda)^{-\frac{k}{2}} \\ & \ -  \frac{d}{d \la} \bigg\{
\sum_{\xi\in\mathscr{P}} c_{\xi}\, (-\la)^{- \xi} \bigg\} +
\left\{ (-\la)^{-1} \sum_{k = 1}^\infty \frac{2^{k-1}
\beta_k}{(\log (-\la))^{k}} \right\} ,
\end{align*}
where the $a_k$ coefficients are independent of $L$, the
$c_{\xi}$'s are the coefficients in \eqref{cxi} and the
$\beta_k$'s are the coefficients in \eqref{p0'p0}.
\end{proposition}
\begin{proof} Since $L = L_0 \oplus L_1$ is decomposable, it follows that
\[
(\Ll_L - \lambda)^{-1} = (\Ll_{L_0} - \lambda)^{-1}  + (\Ll_{L_1}
- \lambda)^{-1}
\]
where $\Ll_{L_0}$ is the operator $\Ll_L$ restricted to the
$-\frac14$ eigenspaces of $A$ and $\Ll_{L_1}$ is the operator
$\Ll_L$ restricted to the eigenspaces of $A$ in $(-\frac14,
\frac34)$.  For the operator $\Ll_{L_1}$,  recall equation
\eqref{e:decom},
\begin{equation*}
\frac{d}{d x} \log F(i x) \, \sim \, \sum_{k=0}^\infty b_k x^{-k}
+ \sum - 2 \xi c_{ \xi}\, x^{-2 \xi - 1}.
\end{equation*}
Combining this with Theorem \ref{thm-res}, we have
\[
\Tr (\Ll_{L_1} - \lambda)^{-1} \sim \sum_{k=1}^\infty d_k \, (-
\lambda)^{-\frac{k}{2}} - \frac{d}{d \la} \bigg\{ \sum c_{\xi}\,
(-\la)^{- \xi} \bigg\}
\]
where the $d_k$ coefficients are independent of $L_1$ and the
$c_{\xi}$'s are the coefficients in \eqref{cxi} or
\eqref{coeffcxi}. Let $F_0(\mu)$ denote the function $F(\mu)$ in
Proposition \ref{prop-Fmu} in the case that $q_1 = 0$, $\Aa =
\Aa_0$, and $\Bb = \Bb_0$ where $\begin{pmatrix} \Aa_0 & \Bb_0
\end{pmatrix}$ defines $L_0$. Then just as in the proof of
Proposition \ref{prop-asympresol}, in conjunction with Proposition
\ref{prop-asympF0} (see \eqref{asymFlog}):
\[
\frac{d}{d x} \log F_0(i x) \, \sim \, \sum_{k = 1}^\infty
\frac{\beta_k}{x(\log x)^k} + \sum_{k = 0}^\infty e_k x^{-k} ,
\]
where the $\beta_k$'s are the coefficients in \eqref{p0'p0}, using
again Theorem \ref{thm-res}, we obtain
\[
\Tr (\Ll_{L_0} - \lambda)^{-1} \, \sim \, \sum_{k=1}^\infty f_k \,
(- \lambda)^{-\frac{k}{2}} \ +\ \left\{ (-\la)^{-1} \sum_{k =
1}^\infty \frac{2^{k-1} \beta_k}{(\log (-\la))^{k}} \right\} .
\]
Combining this with $\Tr (\Ll_{L_1} - \lambda)^{-1}$ analyzed just
before completes our proof.
\end{proof}

\subsection{The heat kernel} \label{sec-hktrace}

\begin{figure}
\centering \includegraphics{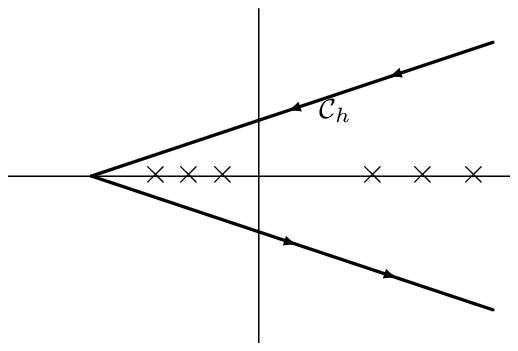} \caption{The
contour $\Cc_h$.} \label{fig-contstrange}
\end{figure}

Now we consider the asymptotics of the trace of $e^{-t\Ll_L}$ as
$t\to 0$. For this, we use
\begin{equation}\label{e:heat-trans}
\Tr (e^{-t \Ll_L} ) = \frac{i}{2 \pi} \int_{\Cc_h} e^{-t \la} \Tr
(\Ll_L - \la)^{-1} d \la
\end{equation}
where $\Cc_h$ is a counter-clockwise contour in the plane
surrounding eigenvalues of $\Ll_L$; see Figure
\ref{fig-contstrange}.
Then the small-time asymptotics of the heat trace is determined by
the large-spectral parameter asymptotics of $\Tr (\Ll_L -
\la)^{-1}$ as we will see in the following proposition.

\begin{proposition} \label{prop-heatk}
For an \textbf{arbitrary} Lagrangian subspace $L \subset \C^{2q}$,
as $t \to 0$ we have
\begin{align*}
\Tr ( e^{-t \Ll_{L}})\, \sim \, & \sum_{k=0}^\infty
\widetilde{a}_k\, t^{\frac{-1+k}{2}}
+\sum_{k=0}^\infty \widetilde{b}_k (\log t)^{-1-k}
\\ & +
\sum_{\xi\in \mathscr{\pP}}\sum^{|p_\xi|+1}_{k=0}
\widetilde{c}_{\xi k}\, t^{\xi} (\log t)^{k} + \sum_{\xi\in
\mathscr{\lL}}\sum_{k=0}^\infty \widetilde{d}_{\xi k}\, t^{\xi}\,
(\log t)^{-\ell_\xi-k} ,
\end{align*}
with $\widetilde c _{10} = 0$ and $\widetilde c_{\xi (|p_\xi| +1)}
=0$ for $\xi \notin \N_0$.
\end{proposition}

\begin{proof}
By Proposition \ref{prop-asympresol},  we have
\begin{multline}\label{e:asymp}
\Tr( \Ll_L - \la)^{-1}\,  \sim \, \sum_{k = 1}^\infty a_k
(-\la)^{-\frac k2}
+ \frac{q_0 - j_0}{(-\la) (\log (-\la) - 2 \widetilde{\gamma})}  \\
-\frac{d}{d \la} \left\{ \sum 2^\ell c_{\ell \xi}\, (-\la)^{- \xi}
\Big(2 \widetilde{\gamma} - \log (-\la) \Big)^{-\ell} \right\}
\end{multline}
as $|\la| \to \infty$ with $\la$ in a sector not intersecting the
positive real axis. We use \eqref{e:heat-trans} for each term on
the right hand side. For the first term, making the change of
variables $\la \mapsto t^{-1}\la$,
\[
\int_{\Cc_h} e^{-t\la}(-\la)^{-\frac k2} d\la = t^{\frac k2 -1}
\int_{t \Cc_h} e^{-\la}(-\la)^{-\frac k2} d\la.
\]
The integral part depends on $t$ via $t\Cc_h$, but is smooth at
$t=0$. Hence, the first part $\sum_{k = 1}^\infty a_k
(-\la)^{-\frac k2}$ contributes
\begin{equation}\label{e:cont1}
\sum_{k=0}^\infty a'_{k}\, t^{\frac{-1+k}{2}}.
\end{equation}
For the third term on the right hand side of \eqref{e:asymp},
using integration by parts, we have
\begin{multline*}
\int_{\Cc_h} e^{-t\la} \frac{d}{d \la} \left\{  (-\la)^{- \xi}
\Big(2 \widetilde{\gamma} - \log
(-\la) \Big)^{-\ell} \right\} \, d\la\\
=t \int_{\Cc_h} e^{-t\la} \left\{ (-\la)^{- \xi} \Big(2
\widetilde{\gamma} - \log (-\la) \Big)^{-\ell} \right\} \, d\la.
\end{multline*}
Deforming $\Cc_h$ to the real line, we find
\begin{eqnarray}\label{e:t}
 & & t \int_{\Cc_h} e^{-t\la} \left\{  (-\la)^{- \xi} \Big(2
\widetilde{\gamma} - \log (-\la)
\Big)^{-\ell} \right\} \, d\la\nonumber\\
& & \hspace{2.0cm} = t \Big( \int^1_\infty e^{-tx} (-(x+i 0))^{-
\xi} \Big(2
\widetilde{\gamma} - \log (-(x+i 0)) \Big)^{-\ell}\, dx\nonumber\\
& & \hspace{3.0cm} +\int^\infty_1 e^{-tx} (-(x-i 0))^{- \xi}
\Big(2 \widetilde{\gamma} -
\log (-(x-i 0)) \Big)^{-\ell}\, dx +h(t) \Big)\nonumber\\
& & \hspace{2.0cm} = t \Big( e^{-i \pi \xi} \int_1^\infty e^{-tx}
x^{- \xi} \Big(2
\widetilde{\gamma} - \log x - i \pi  \Big)^{-\ell}\, dx\nonumber\\
& & \hspace{3.0cm} - e^{i \pi \xi } \int^\infty_1 e^{-tx} x^{-
\xi} \Big(2 \widetilde{\gamma} - \log x +i\pi  \Big)^{-\ell}\, dx
+h(t) \Big) \, ,\nonumber
\end{eqnarray}
where $h(t)$ is a smooth function at $t=0$. Since for any complex
number $z$ we have $i (z-\bar z) = -2 \Im z$, we see that modulo a
term that is a smooth function of $t$ at $t=0$,
\begin{eqnarray}
\frac i {2 \pi} \int_{\Cc_h} e^{-t\la} \frac{d}{d \la} \left\{
(-\la)^{- \xi} \Big(2 \widetilde{\gamma} - \log (-\la)
\Big)^{-\ell} \right\} \, d\la  = - \frac t \pi \Im \ell (t) \,
,\nonumber\end{eqnarray} where
\begin{eqnarray}
 \ell (t) = e^{- i \pi \xi} \int\limits_1^\infty e^{-tx } x^{-\xi}
\left( 2\widetilde \gamma - \log x - i \pi \right) ^{-\ell} dx
;\label{deffell}
\end{eqnarray}
we shall compute the asymptotics of $\ell (t)$ as $t\to 0$. To do
so, let $j\geq \xi > j-1$, $j\in \N_0$; observe that the $j$-th
derivative $\ell ^{(j )} (t)$ of $\ell (t)$ is given by $$ \ell
^{(j)} (t) = e^{- i \pi \xi} (-1)^j \int\limits_1^\infty e^{-tx}
x^{j-\xi}  \left( 2 \widetilde \gamma - \log x -i \pi \right)
^{-\ell} dx .$$ Note that $x^{j-\xi}\cdot(2\widetilde \gamma -\log
x -i \pi )^{-\ell}$ is integrable near $x=0$, so we can write
$$  \ell ^{(j)} (t) = e^{-i \pi \xi} (-1)^j (f(t) + g(t)) $$ with $$f(t) :=
\int\limits_0^\infty e^{-tx} \,\,x^{j-\xi} \left( 2 \widetilde
\gamma - \log x -i \pi \right) ^{-\ell} dx $$ and $$g(t) = -
\int\limits_0^1 e^{-tx} \,\,x^{j-\xi} \left( 2 \widetilde \gamma -
\log x -i \pi \right) ^{-\ell} dx .$$ Note that $g(t)$ is smooth
at $t=0$. We will now determine the asymptotics of $f(t)$ near
$t=0$. To this end, we make the change of variables $x \mapsto
t^{-1} x$:
$$f(t) = t^{\xi -j -1}
\int\limits_0^\infty e^{-x} \,\, x^{j-\xi} \left( 2 \widetilde
\gamma - \log x + \log t -i \pi \right) ^{-\ell} dx \,\, .$$ We
need to consider two cases; $\ell \leq 0$ and $\ell>0$. For $\ell
\leq 0$ we use the binomial expansion to find
\begin{eqnarray}
f(t) &=& t^{\xi -j -1} \sum^{|\ell|}_{k=0} \binom{|\ell|}{k} (\log
t)^{|\ell|-k} \int^\infty_0 e^{-x} \,\, x^{j-\xi}
\left( 2 \widetilde \gamma - \log x - i\pi \right)^k dx \nonumber\\
&=& t^{\xi-j -1} \sum_{k=0}^{|\ell |} c_{\xi, k,\ell} (\log t)^k ,
\nonumber
\end{eqnarray}
with suitable coefficients $c_{k,\ell}$. For $\ell >0$ we first
write $$f(t) = t^{\xi - j -1} (\log t ) ^{-\ell}
\int\limits_0^\infty e^{-x} \,\, x^{j-\xi} \left( 1-\frac{ \log x
+ i\pi - 2 \widetilde \gamma } {\log t} \right)^{-\ell} dx.$$
Since $(1 - r)^{-1} = \sum_{k = 0}^N r^k + r^{N+1} (1 - r)^{-1}$
for any $N \in \N$, we see that for any $N \in \N$,
$$ (1-r)^{-\ell } = \sum_{k=0} ^{N-\ell +1} a_{k,\ell } \,\,r^k
+ r^{N-\ell +2} \,\,\sum_{k=0} ^{\ell -1} b_{k,\ell} \,\,\frac{
r^k} { (1-r) ^{k+1}} .$$ For $f(t)$ this implies
\begin{eqnarray}
f(t) &=& t^{\xi - j -1} (\log t) ^{-\ell} \sum_{k=0}^{N-\ell +1}
a_{k,\ell} (\log t)^{-k} \int\limits_0^\infty
e^{-x}\,\, x^{j-\xi}
\left( \log x + i\pi -2 \widetilde \gamma \right)^k  dx+ \nonumber\\
& &\hspace{-1.5cm}t^{\xi -j -1} (\log t)^{-\ell} \,\, (\log
t)^{-(N-\ell +2)} \sum_{k=0}^{\ell -1} b_{k,\ell } (\log t ) ^{-k}
\int\limits_0^\infty e^{-x} \,\, x^{j-\xi} \frac{ (\log x + i\pi
-2 \widetilde \gamma )^{N-\ell +2+k}} {\left( 1-\frac{\log x +
i\pi -2 \widetilde \gamma}{\log t} \right)^{k+1} } dx
.\nonumber\end{eqnarray} The last integral is bounded as $t\to 0$;
as $N \in \N$ was arbitrary, we conclude for $\ell >0$
\begin{eqnarray} f(t) \sim t^{\xi -j -1} \sum_{k=0}^\infty \widetilde a_{\xi ,k,\ell}
(\log t)^{-k-\ell} .\nonumber\end{eqnarray} In summary: for $\ell
\leq 0$ we have shown
\begin{eqnarray} \ell ^{(j )} (t) \sim e^{-i\pi \xi} \left( t^{\xi - j-1} \sum_{k=0}^{|\ell |}
c_{\xi ,k,\ell} (\log t)^k + \sum _{k=0} ^\infty \gamma _{\xi ,
k,\ell}\,\, t^k  \right), \label{asympls0}\end{eqnarray} while for
$\ell
>0$, we found
\begin{eqnarray} \ell ^{(j )} (t) \sim e^{-i \pi \xi} \left( t^{\xi - j-1} \sum_{k=0}^\infty
\widetilde a _{\xi ,k,\ell} (\log t)^{-k-\ell} + \sum _{k=0}
^\infty \widetilde \gamma _{\xi ,k,\ell} \,\,t^k  \right) .
\label{asymplg0}\end{eqnarray} In order to find the small-$t$
asymptotics of $\ell (t)$ we need to integrate $j$ times. Using
\cite{BGrI-RyI00}, equation 2.722,
\begin{eqnarray}
\int t^n (\log t)^m dt &=& \frac {t^{n+1}}{m+1} \sum_{k=0}^m (-1)
^k \, (m+1) \cdot m \cdot \cdot \cdot (m-k+1) \frac { (\log
t)^{m-k}}{(n+1)^{k+1}},\nonumber \end{eqnarray} for $n\neq -1$,
$m\neq -1$, and \cite{BGrI-RyI00}, equation 2.724, in the form
\begin{eqnarray}
 \int \frac{t^n}{(\log t)^m } dt &=&  \frac{ t^{n+1}}{(n+1) (\log t)^{m}} +
\frac {m}{n+1} \int \frac{t^n}{(\log t)^{m+1} }dt , \nonumber
\end{eqnarray}
for $n\neq -1$, $m\neq 0$, in addition
\[
\int t^{-1} (\log t)^{-1} \, d t = \log | \log t|\quad , \quad
\int t^{-1} (\log t)^{-k - 1} d t = - \frac{1}{k} (\log t)^{-k} \
\ \text{for} \ \ k \neq 0,
\]
we obtain for $\ell \leq 0$, $\xi \notin \N_0$, that
\begin{eqnarray}
\Im \ell (t) \sim t^{\xi -1} \sum_{k=0}^{|\ell | } c'
_{\xi,k,\ell} (\log t)^k + \sum_{k=0} ^\infty \gamma '
_{\xi,k,\ell} t^k , \label{imls0}\end{eqnarray} whereas for $\xi
\in \N_0$ the first summation extends up to $|l|+1$ and
$c'_{1,0,\ell}=0$.

For $l>0$ the answer reads
\begin{eqnarray}
\Im \ell (t) \sim t^{\xi -1} \sum_{k=0}^{\infty} \widetilde c
'_{\xi ,k,\ell} (\log t)^{-k-\ell} + \sum_{k=0} ^\infty \widetilde
\gamma ' _{\xi,k,\ell} t^k  .\label{imlg0}\end{eqnarray}

Contributions from the second term in (\ref{e:asymp}) are found
from (\ref{imlg0}) with $\xi =1$ and $\ell =1$,
\begin{eqnarray} \sum_{k=0}^\infty \widetilde c_k (\log t)^{-k-1}
+ \sum_{k=0} ^\infty \widetilde \gamma _k t^k . \label{extra}
\end{eqnarray}
Combining (\ref{e:cont1}), (\ref{imls0}), (\ref{imlg0}) and
(\ref{extra}) completes the proof.
\end{proof}

Using Proposition \ref{prop-restracede} and repeating the proof of
Proposition \ref{prop-heatk}, we have
\begin{proposition} \label{prop-heatk2}
For an arbitrary decomposable Lagrangian $L$, the heat kernel
$e^{-t \Ll_L}$ has the following trace expansion as $t \to 0$:
\[
\Tr ( e^{-t \Ll_{L}})\, \sim \, \sum_{k=0}^\infty
\widetilde{a}_{k}\, t^{\frac{-1+k}{2}} +
\sum_{\xi\in\mathscr{\pP}} \widetilde{c}_\xi \, t^{\xi} + \sum_{k
= 1}^\infty \widetilde{d}_k(\log t)^{-k}.
\]
\end{proposition}

\section{Proofs of the main theorems}\label{s:meromor}

We now prove our main results starting with the resolvent
expansion.

\subsection{The resolvent expansion --- Theorems \ref{thm-restracegen}
and \ref{thm-restrace}}

We work under the assumptions of Theorem \ref{thm-restracegen}, so
$\Lambda \subset \C$ denotes a sector not intersecting the
positive real axis and $L \subset V$ denotes a given, but
arbitrary, Lagrangian subspace of $V$. We cut the manifold $M$ at
the hypersurface $r = R$ in the collar $[0,\varepsilon)_r \times
\Gamma$ with $0 < R < \varepsilon$, giving a decomposition
\[
M = X \cup Y,
\]
where $X = [0,R]_r \times \Gamma$ and $Y$ is a manifold with a
collared neighborhood $[R,\varepsilon)_r \times \Gamma$ near its
boundary, which we identify with $\Gamma$. Let $\Delta_Y$ denote
the restriction of $\Delta$ to $Y$ with the Dirichlet condition at
$r = R$ and let $\Delta_{X,L}$ denote the restriction of $\Delta$
to $X$:
\[
\Delta_{X,L} : = - \partial_r^2 + \frac{1}{r^2} A_\Gamma
\]
with domain the restriction of $\mathrm{dom}(\Delta_L)$ to $X$ and
with the Dirichlet condition at $r = R$. It is well-known that the
Schwartz kernel of the resolvent $(\Delta_Y - \la)^{-1}(y,y')$,
where $(y,y') \in Y \times Y$, is a \emph{smooth} function of
$(y,y')$ and vanishes to infinite order as $|\la| \to \infty$ with
$\la \in \Lambda$ as long as $y \ne y'$ (see for instance
\cite{SeR69}). In the following lemma we prove a similar statement
for the operator $\Delta_{X,L}$ on the generalized cone.

\begin{lemma} \label{lem-smoothing}
If $\varphi,\psi \in C^\infty (X)$ have disjoint supports, then
for any differential operator $P$ that vanishes near $\partial M$,
the operator
\[
\varphi P (\Delta_{X,L} - \la)^{-1} \psi
\]
is a trace-class operator that vanishes, with all derivatives, to
infinite order (in the trace-class norm) as $|\la| \to \infty$
with $\la \in \Lambda$.
\end{lemma}
\begin{proof}
If we prove this theorem for $\overline{\psi} (\Delta_{X,L} -
\overline{\la})^{-1} P^* \overline{\varphi}$, then taking adjoints
we get our theorem. Hence, we just have to prove the corresponding
statement for $\varphi (\Delta_{X,L} - \la)^{-1} P \psi$, where $P
\psi$ is the operator $f \mapsto P(\psi f)$. We prove this lemma
using the heat kernel $e^{-t \Delta_{X,L}}$, whose structure is
found in \cite{MoE99}. To this end, observe that
\[
(\Delta_{X,L} - \la)^{-1} = \int_0^1 e^{t \la} \, e^{-t
\Delta_{X,L}}\, dt + e^{\la} (\Delta_{X,L} - \la)^{-1}
e^{-\Delta_{X,L}}.
\]
Then
\begin{equation} \label{resXLsm}
\varphi (\Delta_{X,L} - \la)^{-1} P \psi = \int_0^1 e^{t \la} \,
\varphi e^{-t \Delta_{X,L}}P \psi\, dt + e^{\la} \varphi
(\Delta_{X,L} - \la)^{-1} e^{-\Delta_{X,L}} P \psi.
\end{equation}
Assume for the moment that $\Lambda \subset \C$ is contained
entirely in the left-half plane (so that $\Re \la \to -\infty$ as
$|\lambda| \to \infty$ with $\la \in \Lambda$). Now, the operator
$e^{-\Delta_{X,L}} P$ is of trace-class and $(\Delta_{X,L} -
\la)^{-1}$ is a bounded operator which decays like $|\la|^{-1}$ as
$|\la| \to \infty$ with $\la \in \Lambda$. Therefore, since the
trace-class operators form an ideal within the bounded operators,
the operator
\[
(\Delta_{X,L} - \la)^{-1} e^{-\Delta_{X,L}} P
\]
is of trace-class and it decays like $|\la|^{-1}$ as $|\la| \to
\infty$ with $\la \in \Lambda$. Hence,
\[
e^{\la} (\Delta_{X,L} - \la)^{-1} e^{-\Delta_{X,L}} P
\]
decays exponentially, with all derivatives, in the trace-class
operators as $|\la| \to \infty$ with $\la \in \Lambda$ (recall
that $\Re \la \to -\infty$ as $|\lambda| \to \infty$ with $\la \in
\Lambda$). Therefore, the second operator in \eqref{resXLsm}
decays exponentially in the trace-class operators as $|\la| \to
\infty$ with $\la \in \Lambda$. By the main theorem of
\cite{MoE99} (see also Theorem 4.1 of loc.\ cit.),\ since the
supports of $\varphi$ and $P \psi$ are disjoint, it follows that
the operator
\[
\varphi e^{-t \Delta_{X,L}}P \psi
\]
is a trace-class operator that vanishes to infinite order at $t =
0$ (within the trace-class operators). Therefore, the operator
$\int_0^1 e^{t \la} \, \varphi e^{-t \Delta_{X,L}}P \psi\, dt$ in
\eqref{resXLsm} decays exponentially, with all derivatives, in the
trace-class operators as $|\la| \to \infty$ with $\la \in
\Lambda$.

Summarizing: We have proved our theorem when $\Lambda \subset \C$
is contained entirely in the left-half plane (so that $\Re \la \to
-\infty$ as $|\lambda| \to \infty$ with $\la \in \Lambda$). Our
proof is finished once we finish the cases when $\Lambda \subset
\C$ is contained entirely in the upper-half plane and lower-half
plane; for concreteness, let us focus on the upper-half plane.
Then we can fix a complex number $a \in \C$ with \emph{positive}
real part (and positive imaginary part) such that $a \cdot \Lambda
\subset \C$ is entirely contained in the left-half plane. Then one
can construct the heat kernel $e^{-t a \Delta_{X,L}}$ (cf.\
\cite[p.\ 282--284]{BMeR93}) which has the same trace-class
properties as $e^{-t \Delta_{X,L}}$ as described in \cite[Th.\
4.1]{MoE99}. Now we proceed as above: Just as we wrote
\eqref{resXLsm}, one can check that
\begin{equation} \label{resXLsm2}
\varphi (\Delta_{X,L} - \la)^{-1} P \psi = a \int_0^1 e^{t \la} \,
\varphi e^{-t a \Delta_{X,L}}P \psi\, dt + e^{a \la} \varphi
(\Delta_{X,L} - \la)^{-1} e^{-a \Delta_{X,L}} P \psi.
\end{equation}
By the choice of $a$, note that as $|\la| \to \infty$ with $\la
\in \Lambda$, we have $\Re (a \la) \to -\infty$ as $|\lambda| \to
\infty$ with $\la \in \Lambda$. Therefore, analyzing
\eqref{resXLsm2} by repeating the argument we used in the previous
paragraph to analyze \eqref{resXLsm} proves our lemma in the case
when $\Lambda \subset \C$ is contained entirely in the upper-half
plane.
\end{proof}

\begin{figure} \centering
\includegraphics{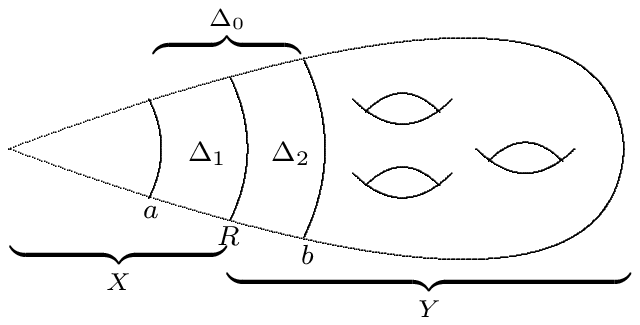} \caption{The maps $\Delta_0,\Delta_1,
\Delta_2$.} \label{fig-conecut1}
\end{figure}

Let us \emph{fix} $0 < a < R$ and $R < b < \varepsilon$, and
define
\[
M_0: = [a,b] \times \Gamma,\quad M_1 := [a,R] \times \Gamma ,\quad
M_2 := [R, b] \times \Gamma.
\]
For $j = 0,1,2$, let $\Delta_j$ denote the Laplacian on $M_j$ with
the Dirichlet boundary condition at the boundaries of $M_j$; see
Figure \ref{fig-conecut1}.
The importance of the operators $\Delta_0,\Delta_1,\Delta_2$ is
that they are \emph{smooth} (not singular) Laplace-type operators
on compact manifolds with boundary with local boundary conditions,
the properties of which are completely understood
\cite{Se66,SeR69}. The idea to prove Theorem \ref{thm-restracegen}
is to compare the resolvents on $M$, $X$, and $Y$ to those on
$M_0$, $M_1$, and $M_2$.

\begin{lemma} \label{lem-differences}
The differences of resolvents
\begin{multline*}
\Ss(\la) := (\Delta_L - \la)^{-1} - (\Delta_{X,L} - \la)^{-1} -
(\Delta_Y - \la)^{-1} -\\ \Big( (\Delta_0 - \la)^{-1} - (\Delta_1
- \la)^{-1} - (\Delta_2 - \la)^{-1} \Big)
\end{multline*}
is trace-class and vanishes, with all derivatives, to infinite
order (in the trace-class norm) as $|\la| \to \infty$ with $\la
\in \Lambda$.
\end{lemma}
\begin{proof} Let $\varrho(r) \in C^\infty(\R)$ be a
non-decreasing function such that $\varrho(r) = 0$ for $r \leq
1/4$ and $\varrho(r) = 1$ for $r \geq 3/4$. For real numbers
$\alpha < \beta$, we define $\varrho_{\alpha,\beta}(r) :=
\varrho\Big( \frac{r-\alpha}{\beta-\alpha} \Big)$. The main
properties of $\varrho_{\alpha,\beta}$ we will use below are that
$\varrho_{\alpha,\beta}(r) = 0$ on a neighborhood of $\{r \leq
\alpha\}$ and $\varrho_{\alpha,\beta}(r) = 1$ on a neighborhood of
$\{r \geq \beta\}$. Let us choose real numbers $a_1,a_2,b_1,b_2$
such that
\[
a < a_1 < a_2 < R < b_1 < b_2 < b.
\]
We define
\[
\psi_1(r) := 1 - \varrho_{a_1,a_2}(r)\, ,\ \ \psi_2(r) :=
\varrho_{b_1,b_2}(r) \, ,\ \ \psi_0(r) := 1 - \psi_1(r) -
\psi_2(r) ,
\]
and
\[
\varphi_1(r) := 1  - \varrho_{a_2,R}(r) \, ,\ \ \varphi_2(r) :=
\varrho_{R,b_1}(r) \, ,\ \ \varphi_0(r) := 1 - \varrho_{a,a_1}(r)
- \varrho_{b_2,b}(r).
\]
The functions $\{\psi_i\}$, $\{\varphi_i\}$ extend either by $0$
or $1$ to define smooth functions on all of $M$ and $\{\psi_i\}$
forms a partition of unity of $M$ such that $\varphi_i = 1$ on
$\mathrm{supp}(\psi_i)$. Now to prove this lemma, we first claim
that each of the following equalities holds modulo a trace-class
operator vanishing to infinite order  as $|\la| \to \infty$ with
$\la \in \Lambda$:
\begin{equation}
\label{sumres}
\begin{split}
(\Delta_L - \la)^{-1} & = \varphi_1 (\Delta_{X,L} - \la)^{-1}
\psi_1 + \varphi_0 (\Delta_0 - \la)^{-1}
\psi_0 + \varphi_2 (\Delta_Y - \la)^{-1} \psi_2 ,\\
(\Delta_{X,L} - \la)^{-1} & = \varphi_1 (\Delta_{X,L} - \la)^{-1}
\psi_1 + \varphi_0 (\Delta_1 - \la)^{-1}\psi_0 ,\\
(\Delta_Y - \la)^{-1} & = \varphi_0 (\Delta_2 - \la)^{-1}\psi_0 +
\varphi_2 (\Delta_Y - \la)^{-1} \psi_2 .
\end{split}
\end{equation}
For instance, let us verify the first claim in \eqref{sumres}; the
other claims are verified using a similar argument. Define
\[
Q(\la) := \varphi_1 (\Delta_{X,L} - \la)^{-1} \psi_1 + \varphi_0
(\Delta_0 - \la)^{-1}\psi_0 + \varphi_2 (\Delta_Y -
\la)^{-1}\psi_2.
\]
Then observe that $(\Delta_L - \la) Q(\la) = \Id + K(\la)$ where
\begin{multline*}
K(\la) =
[ (\Delta_L  - \la), \varphi_1] (\Delta_{X,L} - \la)^{-1} \psi_1 \\
+ [ (\Delta_L - \la), \varphi_0] (\Delta_0 - \la)^{-1} \psi_0 + [
(\Delta_L - \la), \varphi_2] (\Delta_Y - \la)^{-1} \psi_2,
\end{multline*}
where $[\ ,\ ]$ denotes the ``commutator". Now, because the
supports of
\[
[ (\Delta_L - \la), \varphi_i] = - \left[ \partial_r^2 , \varphi_i
\right] = -\big(\varphi_i'' + 2 \varphi'\, \partial_r\big)
\]
and $\psi_i$ are disjoint, it follows that each of the three
operators making up $K(\la)$ is trace-class and vanishes to
infinite order (in the trace-class norm) as $|\la| \to \infty$
with $\la \in \Lambda$; indeed, this statement for $[ (\Delta_L  -
\la), \varphi_1] (\Delta_{X,L} - \la)^{-1} \psi_1$ follows from
Lemma \ref{lem-smoothing} and the statements for $[ (\Delta_L -
\la), \varphi_0] (\Delta_0 - \la)^{-1} \psi_0$ and $[ (\Delta_L -
\la), \varphi_2] (\Delta_Y - \la)^{-1} \psi_2$ are well-known (see
e.g.\ \cite{SeR69}). Therefore, $K(\la)$ is trace-class and
vanishes to infinite order (in the trace-class norm) as $|\la| \to
\infty$ with $\la \in \Lambda$. Now applying $(\Delta_L -
\la)^{-1}$ to both sides of $(\Delta_L - \la) Q(\la) = \Id +
K(\la)$, we obtain
\[
(\Delta_L - \la)^{-1} =  Q(\la) - (\Delta_L - \la)^{-1} K(\la),
\]
which establishes our claim for the first equality in
\eqref{sumres}. A similar argument works to prove that the other
equalities in \eqref{sumres} hold modulo trace-class with infinite
decay (with all derivatives as $|\la| \to \infty$ with $\la \in
\Lambda$). From \eqref{sumres}, it follows that modulo trace-class
with infinite decay,
\begin{multline} \label{DeltaM12}
(\Delta_L - \la)^{-1} - (\Delta_{X,L} - \la)^{-1} - (\Delta_Y -
\la)^{-1}\\ = \varphi_0 (\Delta_0 - \la)^{-1}\psi_0 - \varphi_0
(\Delta_1 - \la)^{-1}\psi_0 - \varphi_0 (\Delta_2 -
\la)^{-1}\psi_0.
\end{multline}

On the other hand, very similar arguments used to establish
\eqref{sumres} show that modulo trace-class with infinite decay,
\[
\begin{split}
(\Delta_0 - \la)^{-1} & = \varphi_1 (\Delta_1 - \la)^{-1} \psi_1 +
\varphi_0 (\Delta_0 - \la)^{-1}\psi_0 +  \varphi_2 (\Delta_2 -
\la)^{-1} \psi_2, \\ (\Delta_1 - \la)^{-1} & = \varphi_1
(\Delta_1 - \la)^{-1} \psi_1 + \varphi_0 (\Delta_1 - \la)^{-1}\psi_0, \\
(\Delta_2 - \la)^{-1} & = \varphi_0 (\Delta_2 - \la)^{-1}\psi_0 +
\varphi_2 (\Delta_2 - \la)^{-1} \psi_2 .
\end{split}
\]
Combining these identities we can write, modulo trace-class with
infinite decay,
\begin{multline*}
(\Delta_0 - \la)^{-1} - (\Delta_1 - \la)^{-1} - (\Delta_2 - \la)^{-1}\\
= \varphi_0 (\Delta_0 - \la)^{-1}\psi_0 - \varphi_0 (\Delta_1 -
\la)^{-1}\psi_0 - \varphi_0 (\Delta_2 - \la)^{-1}\psi_0.
\end{multline*}
Comparing this with \eqref{DeltaM12} completes the proof of our
lemma.
\end{proof}

Using our standard notation, let $\{\lambda_\ell\}$ denote the set
of all eigenvalues of $A_\Gamma$ and let $E_\ell$ denote the span
of the $\la_\ell$-th eigenvector. Let $\Pi$ and $\Pi^\perp$
denote, respectively, the orthogonal projections of
$L^2(\Gamma,E_\Gamma)$ onto $W : = \bigoplus_{-\frac14 \leq
\lambda_\ell < \frac34} E_\ell$ and $W^\perp$. Using the isometry
between
\[
L^2([0,R]\times \Gamma, E) \cong L^2([0,R], L^2(\Gamma, E_\Gamma))
\]
where $E_\Gamma:=E|_\Gamma$, we obtain the corresponding
projections on $L^2([0,R]\times \Gamma, E)$, which we denote by
the same notations $\Pi$ and $\Pi^\perp$. Let  $\Ll_L$ denote the
model operator introduced in Section \ref{s:modelproblems}
(specifically, Section \ref{subs:model}), and define
\[
\Delta'_{X} := - \partial_r^2 + \frac{1}{r^{2}} A'_\Gamma,\quad
\text{where\ } A'_\Gamma :=
\begin{cases} \frac{3}{4} & \text{over}\ W\\
A_\Gamma & \text{over}\ W^\perp \end{cases} .
\]

\begin{proposition} \label{prop-import} We have
\begin{align*}
(\Delta_L - \la)^{-1} & =  \Pi (\Ll_L - \la)^{- 1} \Pi + \Pi^\perp
(\Delta_{X}' - \la)^{- 1} \Pi^\perp + (\Delta_Y - \la)^{-1} +\\ &
\hspace{6em} (\Delta_0 - \la)^{-1} - (\Delta_1 - \la)^{-1} -
(\Delta_2 - \la)^{-1} + \Ss(\la)
\end{align*}
where $\Ss(\la)$ is trace-class and vanishes, with all
derivatives, to infinite order as $|\la| \to \infty$ with $\la \in
\Lambda$.
\end{proposition}
\begin{proof}
Observe that
\[
(\Delta_{X,L} - \la)^{- 1} = \Pi (\Delta_{X,L} - \la)^{- 1} \Pi +
\Pi^\perp (\Delta_{X,L} - \la)^{- 1} \Pi^\perp ,
\]
since $\Delta_{X,L}$ preserves $W$ and $W^\perp$, and
\[
\Pi (\Delta_{X,L} - \la)^{- 1} \Pi = \Pi (\Ll_L - \la)^{- 1} \Pi .
\]
Also observe that
\[
\Pi^\perp (\Delta_{X,L} - \la)^{- 1} \Pi^\perp = \Pi^\perp
(\Delta_{X}' - \la)^{- 1} \Pi^\perp.
\]
Hence,
\[
(\Delta_{X,L} - \la)^{- 1} = \Pi (\Ll_L - \la)^{- 1} \Pi +
\Pi^\perp (\Delta_{X}' - \la)^{- 1} \Pi^\perp .
\]
Now solving for $(\Delta_L - \la)^{-1}$ in Lemma
\ref{lem-differences}, we obtain
\begin{align*}
(\Delta_L - \la)^{-1} & =  (\Delta_{X,L} - \la)^{-1} + (\Delta_Y -
\la)^{-1} +\\ & \hspace{6em} (\Delta_0 - \la)^{-1} - (\Delta_1 -
\la)^{-1} - (\Delta_2 - \la)^{-1} + \Ss(\la) \\
& =  \Pi (\Ll_L - \la)^{- 1} \Pi + \Pi^\perp (\Delta_{X}' -
\la)^{- 1} \Pi^\perp + (\Delta_Y - \la)^{-1} +\\ & \hspace{6em}
(\Delta_0 - \la)^{-1} - (\Delta_1 - \la)^{-1} - (\Delta_2 -
\la)^{-1} + \Ss(\la)
\end{align*}
where $\Ss(\la)$ is trace-class and vanishes, with all
derivatives, to infinite order as $|\la| \to \infty$ with $\la \in
\Lambda$. This completes our proof.
\end{proof}

We can now prove Theorem \ref{thm-restracegen}. Let $N \geq
\frac{n}{2}$ with $n = \dim M$. Then taking $N$ derivatives of
both sides of the preceding equality we see that
\begin{multline} \label{resN1}
(\Delta_L - \la)^{-N - 1} =  \Pi (\Ll_L - \la)^{-N - 1} \Pi +
\Pi^\perp (\Delta_{X}' - \la)^{-N - 1} \Pi^\perp \\ + (\Delta_Y -
\la)^{-N - 1} + (\Delta_0 - \la)^{-N - 1} - (\Delta_1 - \la)^{-N -
1}\\ - (\Delta_2 - \la)^{-N - 1} + \Ss^{(N)}(\la) ,
\end{multline}
where $\Ss^{(N)}(\la) = \frac{d^N}{d \la^N} \Ss(\la)$. We now
analyze each term on the right. First, taking $N$ derivatives in
the asymptotic expression of Proposition \ref{prop-asympresol}, we
know that as $|\la| \to \infty$ with $\la \in \Lambda$ we have
\begin{multline} \label{tr1}
\Tr( \Ll_L - \la)^{-N-1}\,  \sim \, \sum_{k = 0}^\infty b_k
(-\la)^{\frac{1- k}2 - N-1} + \frac{1}{N!}\frac{d^N}{d \la^N}
\bigg\{ \frac{q_0 - j_0}{(-\la) (\log (-\la) - 2
\widetilde{\gamma})} \bigg\} \\ - \frac{1}{N!}\frac{d^{N+1}}{d
\la^{N+1}} \bigg\{ \sum 2^\ell c_{\ell \xi}\, (-\la)^{- \xi}
\Big(2 \widetilde{\gamma} - \log (-\la) \Big)^{-\ell} \bigg\} .
\end{multline}
It follows from \cite{LoPR02} that the operator $\Pi^\perp
(\Delta_{X}' - \la)^{-N - 1} \Pi^\perp$ is trace-class and
\begin{equation} \label{tr2}
\Tr( \Pi^\perp (\Delta_{X}' - \la)^{-N - 1} \Pi^\perp) \, \sim \,
\sum_{k=0}^\infty c_k \, (- \lambda)^{\frac{n-k}{2}  - N-1}  + b
\,  (- \lambda)^{- N - 1} \, \log (- \lambda) ;
\end{equation}
in principle (see \cite[Ch.\ 7]{BMeR93}), one can derive this
resolvent expansion with a lot of work from the corresponding heat
kernel expansion \cite{BS1,BS2,Ch2,Chou}. From the work of Seeley
\cite{SeR69}, we also know that each of $(\Delta_Z - \la)^{-N -
1}$, where $Z = Y, 0, 1, 2$, is trace class, and
\begin{equation} \label{tr3}
\Tr( (\Delta_Z - \la)^{-N - 1} ) \, \sim \, \sum_{k=0}^\infty
c_{Z,k} \, (- \lambda)^{\frac{n-k}{2}  - N-1}.
\end{equation}
Finally, we know that $\Ss^{(N)}(\la)$ is trace-class and
vanishes, with all derivatives, to infinite order as $|\la| \to
\infty$ with $\la \in \Lambda$. In conclusion, in view of the
expression \eqref{resN1} and our discussions around \eqref{tr1},
\eqref{tr2} and \eqref{tr3}, we see that $(\Delta_L - \la)^{-N -
1}$ is trace-class, and
\begin{multline*}
\Tr (\Delta_L - \lambda)^{-N-1} \,  \sim \, \sum_{k=0}^\infty a_k
\, (- \lambda)^{\frac{n-k}{2}  - N-1}  + b \,  (- \lambda)^{- N -
1} \, \log (- \lambda)\\
+ \frac{1}{N!}\frac{d^N}{d \la^N} \bigg\{ \frac{q_0 - j_0}{(-\la)
(\log (-\la) - 2 \widetilde{\gamma})} \bigg\} \\ -\frac{1}{N!}
\frac{d^{N+1}}{d \la^{N+1}} \bigg\{ \sum 2^\ell c_{\ell \xi}\,
(-\la)^{- \xi} \Big(2 \widetilde{\gamma} - \log (-\la)
\Big)^{-\ell} \bigg\} .
\end{multline*}
This completes the proof of Theorem \ref{thm-restracegen}. Theorem
\ref{thm-restrace} is established by replacing the trace expansion
\eqref{tr1} with the trace expansion found in Proposition
\ref{prop-restracede}.

\subsection{Proofs of Theorems
\ref{thm-main}, \ref{main thm}, \ref{main thm2}}

We now prove the $\zeta$-function theorem. We start off with
Proposition \ref{prop-import}, which states that
\begin{align*}
(\Delta_L - \la)^{-1} & =  \Pi (\Ll_L - \la)^{- 1} \Pi + \Pi^\perp
(\Delta_{X}' - \la)^{- 1} \Pi^\perp + (\Delta_Y - \la)^{-1} +\\ &
\hspace{6em} (\Delta_0 - \la)^{-1} - (\Delta_1 - \la)^{-1} -
(\Delta_2 - \la)^{-1} + \Ss(\la)
\end{align*}
where $\Ss(\la)$ is trace-class and vanishes, with all
derivatives, to infinite order as $|\la| \to \infty$ with $\la \in
\Lambda$. Therefore, by the definition of the $\zeta$-function:
\[
\zeta(s,\Delta_L) := \Tr (\Delta_L^{-s}P_0^\perp)\ , \quad
\Tr(\Delta_L^{-s}P_0^\perp) = \sum_{\la_\ell<0} \la_\ell^{-s}+
\int_{\Re \la = \delta} \la^{-s} (\Delta_L - \la)^{-1} d \la,
\]
where $P_0$ is the orthogonal projection onto $\ker \Delta_L$ and
$\delta
> 0$ is any positive number sufficiently small so that the
spectrum of $\Delta_L$ intersected with $(0,\delta]$ is empty, it
follows that
\begin{equation} \label{zetasDeltaL}
\zeta(s,\Delta_L) \equiv \zeta(s, \Ll_L) + \zeta(s, \Delta_X') +
\zeta(s,\Delta_Y) + \zeta(s, \Delta_0) - \zeta(s,\Delta_1) -
\zeta(s,\Delta_2)
\end{equation}
modulo an entire function. The $\zeta$-function $\zeta(s,\Ll_L)$
is studied thoroughly in Proposition \ref{prop-main}. Also, by the
standard relation between the asymptotics of the resolvent and the
poles of the $\zeta$-function (see e.g.\ \cite{GS96}) it follows
from the resolvent expansions \eqref{tr2} and \eqref{tr3} that
$\zeta(s, \Delta_X')$, $\zeta(s,\Delta_Y)$, $\zeta(s, \Delta_0)$,
$\zeta(s,\Delta_1)$, and $\zeta(s,\Delta_2)$ have the ``regular"
poles at the ``usual" locations $s = \frac{n - k}{2}\notin -\N_0$
for $k \in \N_0$ and, only for $\zeta(s, \Delta_X')$, at $s = 0$
if $\dim \Gamma > 0$. These facts together with Proposition
\ref{prop-main} prove Theorem \ref{thm-main}. Note that Theorems
\ref{main thm} and \ref{main thm2} follow from applying
Proposition \ref{prop-zeta} and Corollary \ref{cor-zeta} to
$\zeta(s, \Ll_L)$ in \eqref{zetasDeltaL}.

\subsection{Proof of Theorems \ref{t:traceM} and \ref{t:traceMde}}

Finally, it remains to prove the heat expansion. As with the proof
for the $\zeta$-function, we start off with Proposition
\ref{prop-import}:
\begin{align*}
(\Delta_L - \la)^{-1} & =  \Pi (\Ll_L - \la)^{- 1} \Pi + \Pi^\perp
(\Delta_{X}' - \la)^{- 1} \Pi^\perp + (\Delta_Y - \la)^{-1} +\\ &
\hspace{6em} (\Delta_0 - \la)^{-1} - (\Delta_1 - \la)^{-1} -
(\Delta_2 - \la)^{-1} + \Ss(\la)
\end{align*}
where $\Ss(\la)$ is trace-class and vanishes, with all
derivatives, to infinite order as $|\la| \to \infty$ with $\la \in
\Lambda$. Then by the definition of the heat operator:
\[
e^{- \Delta_L} := \frac{i}{2 \pi} \int_{\Cc_h} e^{-t \la}
(\Delta_L - \la)^{-1} d \la,
\]
where $\Cc_h$ is a contour as in Figure \ref{fig-contstrange}, it
follows that
\begin{align*}
e^{-t \Delta_L} & =  \Pi e^{-t \Ll_L} \Pi + \Pi^\perp e^{-t
\Delta_{X}'} \Pi^\perp + e^{-t \Delta_Y} +\\ & \hspace{6em} e^{-t
\Delta_0} - e^{-t \Delta_1} - e^{-t \Delta_2} + \Tt(t)
\end{align*}
where $\Tt(t)$ is trace-class and smooth at $t = 0$. Hence,
\begin{multline*}
\Tr( e^{-t \Delta_L} ) = \Tr( e^{-t \Ll_L}) + \Tr( e^{-t
\Delta_X'}) + \Tr(e^{-t\Delta_Y}) \\ + \Tr( e^{-t \Delta_0}) -
\Tr(e^{-t \Delta_1}) - \Tr(e^{-t \Delta_2})
\end{multline*}
modulo a function that is smooth at $t = 0$. The heat trace $\Tr(
e^{-t \Ll_L})$ is studied thoroughly in Proposition
\ref{prop-heatk} and for decomposable Lagrangians in Proposition
\ref{prop-heatk2}. Also, by the standard relation between the
asymptotics of the resolvent and the heat trace expansion (see
e.g.\ \cite{GS96} or Section \ref{sec-hktrace}) it follows from
the resolvent expansions \eqref{tr2} and \eqref{tr3} that $\Tr(
e^{-t \Delta_X'})$, $\Tr(e^{-t\Delta_Y})$, $\Tr( e^{-t
\Delta_0})$, $\Tr(e^{-t \Delta_1})$ and $\Tr(e^{-t \Delta_2})$
have the ``regular" expansion except that $\Tr( e^{-t\Delta_X'})$
may have a $\log t$ term if $\dim \Gamma
> 0$. These facts together with Propositions \ref{prop-heatk}
and \ref{prop-heatk2} prove Theorems \ref{t:traceM} and
\ref{t:traceMde}.

\bibliographystyle{amsplain}
\def\cprime{$'$} \def\polhk#1{\setbox0=\hbox{#1}{\ooalign{\hidewidth
  \lower1.5ex\hbox{`}\hidewidth\crcr\unhbox0}}}
\providecommand{\bysame}{\leavevmode\hbox
to3em{\hrulefill}\thinspace}
\providecommand{\MR}{\relax\ifhmode\unskip\space\fi MR }
\providecommand{\MRhref}[2]{%
  \href{http://www.ams.org/mathscinet-getitem?mr=#1}{#2}
} \providecommand{\href}[2]{#2}

\end{document}